\documentclass[12pt]{amsart}
\let\oldemptyset\emptyset
\def\oldemptyset{o}
\def\oldemptyset{0}

\usepackage[square,compress,comma, numbers,sort]{natbib}
\usepackage[colorlinks=true, citecolor=red, linkcolor=blue]{hyperref}
\usepackage{amsfonts,mathtools}

\allowdisplaybreaks[4]

\def\ThetaH{\Theta_h}
\def\ThetaAH{\Theta_{a+h}}
\def\SZ{ \mathcal{S}_Z}
\def\MZ{ \mathcal{M}_Z}

\usepackage{amssymb}
\usepackage{color}
\newcommand{\abs}[1]{\left\lvert #1 \right\rvert}

\DeclarePairedDelimiterXPP\pk[1]{\mathbb{P}}\{ \}{}{ #1}
\DeclarePairedDelimiterXPP\E[1]{\mathbb{E}}\{ \}{}{	#1}

\def\FRE{\mbox{Fr\'{e}chet }}

\def\toas{\overset{a.s.}\rightarrow}
\def\toprob{\overset{p}\rightarrow}
 \usepackage{xparse}

\NewDocumentCommand{\ceil}{s O{} m}{%
  \IfBooleanTF{#1} 
    {\left\lceil#3\right\rceil} 
    {#2\lceil#3#2\rceil} 
}
\NewDocumentCommand{\floor}{s O{} m}{%
  \IfBooleanTF{#1} 
    {\left\lfloor#3\right\rfloor}
    {#2\lfloor#3#2\rfloor}
}

\def\Hxn{H(x )}
\def\HH{{[h]}}

\newcommand{\norm}[1]{\lVert #1 \rVert}

\definecolor{c20}{rgb}{0.,0.7,0.}
\definecolor{c30}{rgb}{0.,0.,1.}
\definecolor{c40}{rgb}{1,0.1,0.7}
\definecolor{c50}{rgb}{1,0,0}
\definecolor{c60}{rgb}{1,0.9,0.1}
\definecolor{c70}{rgb}{0.50,1.00,0.00}

\def\Ec#1{{\textcolor{c50}{#1}}}
\def\Ec#1{#1}

\def\cE#1{{\textcolor{c30}{#1}}}
\def\cE#1{#1}

\numberwithin{equation}{section}
\newtheorem{theo}{Theorem}[section]
\newtheorem{sat}[theo]{Proposition}
\newtheorem{de}[theo]{Definition}
\newtheorem{lem}[theo]{Lemma}

\newtheorem{example}[theo]{Example}
\newtheorem{korr}[theo]{Corollary}
\newtheorem{remark}[theo]{Remark}

\numberwithin{equation}{section}

\newcommand{\prooftheo}[1]{ \textsc{Proof of Theorem} \ref{#1} }
\newcommand{\proofprop}[1]{\textsc{Proof of Proposition} \ref{#1}}
\newcommand{\prooflem}[1]{\textsc{Proof of Lemma} \ref{#1}}
\newcommand{\proofkorr}[1]{\textsc{Proof of Corollary} \ref{#1}}

\newcommand{\QED}{\hfill $\Box$}

\newcommand{\COM}[1]{}
\newcommand{\COMB}[1]{#1}
\newcommand{\COMC}[1]{#1}

\def\IF{\infty}

\newcommand{\R}{\mathbb{R}}
\newcommand{\inr}{\in \R}

\newcommand{\BQN}{\begin{eqnarray}}
\newcommand{\EQN}{\end{eqnarray}}
\newcommand{\BQNY}{\begin{eqnarray*}}
\newcommand{\EQNY}{\end{eqnarray*}}
\def\ldot{, \ldots,}

\newcommand{\limit}[1]{\lim_{#1 \to   \infty}}

\def\todis{\overset{d}\rightarrow}

\def\bqny#1{ { \begin{eqnarray*} #1 \end{eqnarray*}}}
\def\bqn#1{ { \begin{eqnarray} #1 \end{eqnarray}}}

\newcommand{\BS}{\begin{sat}}
\newcommand{\ES}{\end{sat}}
\newcommand{\BT}{\begin{theo}}
\newcommand{\ET}{\end{theo}}
\newcommand{\BK}{\begin{korr}}
\newcommand{\EK}{\end{korr}}
\newcommand{\EQD}{\stackrel{d}{=}}
\newcommand{\BEX}{\begin{example}}
\newcommand{\EEX}{\end{example}}

\newcommand{\BD}{\begin{de}}
\newcommand{\ED}{\end{de}}

\newcommand{\BIT}{\begin{itemize}}
\newcommand{\EIT}{\end{itemize}}
\newcommand{\BDI}{\begin{description}}
\newcommand{\EDI}{\end{description}}

\newcommand{\BRM}{\begin{remark}}
\newcommand{\ERM}{\end{remark}}

\newcommand{\BEL}{\begin{lem}}
\newcommand{\EEL}{\end{lem}}

\newcommand{\nelem}[1]{{Lemma \ref{#1}}}

\newcommand{\netheo}[1]{{Theorem \ref{#1}}}
\newcommand{\neexamp}[1]{{Example \ref{#1}}}

\newcommand{\equaldis}{\stackrel{{d}}{=}}
\newcommand{\fdd}{\stackrel{\Ec{fdd}}{=}}

\def\TT{\mathcal{T}}
\def\TT{\R }
\def\eZ{\zeta_Z}
\def\Dhz{\Xi_h Z}
\def\thhh{L^h}

\def\HZD{ \mathcal{H}_{Z}^\delta}
\def\Zd{\mathbb{Z}^d}
\def\Zdd{\delta \mathbb{Z}^d}
\def\Z{\mathbb{Z}}
\def\inn{\in \mathbb{N}}

\begin{document}

\title{Representations of max-stable processes via exponential tilting}

\author{Enkelejd  Hashorva}
\address{Enkelejd Hashorva, Department of Actuarial Science 
University of Lausanne,\\
UNIL-Dorigny, 1015 Lausanne, Switzerland}
\email{Enkelejd.Hashorva@unil.ch}

\bigskip

\date{\today}
 \maketitle

\begin{quote}

\COM{
	
    Title: Representations of Max-Stable Processes via Exponential Tilting

    Abstract: The recent contribution Dieker & Mikosch (2015) [1] derived important  representations of max-stable
    Brown-Resnick random fields $\zeta_Z$  with a spectral representation determined by a Gaussian random field  $Z$.
    With motivations from [1] we derive for some general $Z$,
    representations for $\zeta_Z$ via exponential tilting of $Z$ and present further  conditions for the stationarity of $\zeta_Z$.
    Applications of our findings concern  a) Dieker-Mikosch representations  of max-stable processes, b)
    two-sided extensions of stationary max-stable processes, c) tractable representations of max-stable distributions in terms of inf-argmax representation, d) new sufficient conditions for max-domains of attraction, and  e) calculation of generalised Pickands constants. As a by-product we derive a new characterisation of Gaussian random vectors.

    stationary Brown-Resnick random fields with a spectral representation determined by a Gaussian random field  Z.
    With motivations from Dieker & Mikosch (2015) we derive for some general Z    representations for max-stable random fields via exponential tilting of Z and present further new  conditions for the stationarity of those processes.    We shall discuss several applications of our findings which include a) Dieker-Mikosch representations of max-stable processes,    b) two-sided extensions of stationary max-stable processes, c) derivations of tractable representations of max-stable distributions,
    d) sufficient conditions for max-domains of attraction, and  e) calculation of generalised Pickands constants.
    }

{\bf Abstract:} The recent contribution \cite{DM} obtained representations of max-stable stationary Brown-Resnick process $\eZ(t),t\inr^d$  with spectral process $Z$ being  Gaussian. With motivations from \cite{DM} we derive for general $Z$, representations for $\eZ$ via exponential tilting of $Z$. Our  findings concern   {\it Dieker-Mikosch representations} of max-stable processes,   two-sided extensions of stationary max-stable processes, inf-argmax representation of max-stable distribution, and  new formulas for generalised Pickands constants.  Our applications include  conditions for the stationarity of $\eZ$, a characterisation of Gaussian distributions and an alternative proof of Kabluchko's characterisation of Gaussian processes with stationary increments.

\end{quote}


\def\TTT{\mathcal{T}}
 \def\TT{\mathcal{T}}

\section{Introduction}
A random process $\zeta(t),t\in \TT$ is max-stable if all its finite dimensional distributions (fidi's) are max-stable. For simplicity we shall assume hereafter that $\zeta(t)$ has unit Gumbel distribution $e^{-e^{-x}},x\inr$ for any $t\in \TT$ and shall consider $\TT=\R^d$ or $\TT=\Z^d, d\ge 1$. In view of \cite{deHaan} any stochastically continuous max-stable process $\zeta(t),t\in \TT$
satisfies (below $\fdd$ means equality of all fidi's)
$$ \zeta \fdd \eZ,$$
with
 \BQN \label{W}
\eZ(t)&=& \max_{i\ge 1} \Bigl(P_i+ Z_i(t)\Bigr), \quad t\in \TTT,
\EQN
where $Z(t),t\in \TTT $ is  a random process taking values in $[- \IF, \IF)$ with $\E{e^{Z(t)}}=1,t\in \TT$ and 
$\Pi= \sum_{i=1}^\IF \varepsilon_{P_i}$ is a Poisson point process (PPP) on $\R$ with intensity $e^{-x} dx$. Further,  $Z_i$'s are independent copies of $Z$ being also independent of $\Pi$. See for more details the important contributions \cite{Penrose, Schlather2002, MR2324891,kab2009a, KabExt, MolchanovBE,
	Roy, Strokorb, kabStoev,EmbRob,MR3561100}.\\
We shall refer to $\eZ$ as the  {\it associated max-stable} process of $Z$; commonly $Z$ is referred to as the {\it spectral process}. For convenience, we shall write $Z$ as
\bqn{\label{kanon}
	Z(t)= B(t)-  \ln \E*{e^{B(t)}}, \quad t\in \TT,
}
with  $B(t),t \in \mathcal{T} $ a random process  satisfying $\E{e^{B(t)}} < \IF, t \in \TT $. Consequently, $\E{e^{Z(t)}}=1, t\in \TTT$ implying that the marginal distribution functions (df's)  of $\eZ$ are unit Gumbel. \\
One canonical instance is the classical Brown-Resnick construction with $B$ a centred Gaussian process with covariance function $r$ and thus $2 \ln \E{e^{B(t)}}= r(t,t)=:\sigma^2(t),t\in \TTT$.   In view of \cite{kabluchko2011extremes}  the law of $\eZ$ is determined by the incremental variance function $\gamma(s,t)= Var(B(t)- B(s)),s,t\in  \TT $. This fact can be shown by utilising the {\it tilted spectral process}  $\Dhz, h\in \TT$ defined by
\bqny{
 \Dhz(t)&=& B(t)- B(h)- \gamma(h,t)/2 , \quad t\in \TTT.
 }
{\Ec{The law of  $\Dhz$ is uniquely determined by the following conditions: $\Dhz$ is Gaussian,  $\Dhz(h)=0$ almost surely (a.s.)
	and the incremental variance function of $\Dhz$ is $\gamma$. Note that these conditions do
	not involve $\sigma^2$.} \\
\Ec{Next,
setting  $Z^ \HH (t)= B(t)- \sigma^2(t)/2+ r(h,t)$ we have}
\bqn{
 \Dhz(t)&=&  Z^ \HH (t)- \cE{Z^ \HH }(h), \quad  \quad t\in \TT. \label{defDhz}
 }
In view of \nelem{tiltGauss} below $Z ^\HH $ is the exponential tilt of $Z$ by $Z(h)$ i.e.,
\def\pF{\varphi(h)}
$$\pk{Z ^\HH  \in A}= \E*{e^{Z(h)}\mathbb{I}\{Z \in A\}}
, \quad
\forall A \in \mathcal{B}(\R^\mathcal{T}),$$
where $\mathcal{B}(\R^\mathcal{T})$ is the $\sigma$-field generated by all evaluation maps.
 The representation  \eqref{W} implies that  (see e.g., \cite{MolchanovSPA,DM})  
\bqn{
	- \ln  \pk{\eZ(t_i)\le x_i,1 \le i\le n} %
	&=& \E*{ \max_{1 \le i \le n} e^{Z(t_i)- x_i}} \label{tDi}
	= \E*{ e^{Z(h)} \max_{1 \le i \le n} e^{Z(t_i)- Z(h) -x_i}}\notag \\
	&=& \E*{ \max_{1 \le i \le n} e^{\Dhz(t_i) - x_i}}
}
holds for $t_i\in \TTT, x_i \inr, i\le n$ i.e.,
\BQN \label{repA}
\eZ \fdd  \zeta_{\Dhz}.
\EQN
\Ec{Since as mentioned above the process $\Dhz$ can be characterised without making reference to $\sigma^2$, by \eqref{repA} it follows that the law of $\eZ$  depends on $\gamma$ only!}\\
\Ec{Observe that we can define $Z ^\HH $ via exponential tilting for any  process $Z$ such that $\E{e^{Z(h)}}=1$. Furthermore,
		  the calculation of the fidi's of $\eZ$ via \eqref{tDi} does not relate to the Gaussianity of $Z$, but only to the representation \eqref{W} and the fact that}
  \bqn{\label{zh:finite}
	\pk{Z(h)> - \IF }&=&1.
}
 Consequently,  under \eqref{zh:finite} we have that \eqref{repA} is valid for a general {\it spectral process} $Z$ with values in $\R$. Clearly,  \eqref{zh:finite} implies
\bqn{ \label{zho}
\Dhz(h)= 0
}
almost surely, which in view of \cite{HaanPickands}[Lemma 4.1] is a crucial uniqueness condition.\\
The change of measure technique, or in our case the exponential tilting has been utilised in the context of max-stable processes in \cite{kab2009a, MolchanovBE, MR3024101, MR3034327}. In this contribution we present some further developments and applications  that are summarised below:
\BIT
\item[A)] \Ec{According to  \cite{kab2009} the {\it spectral process}} $Z$ is called {\it Brown-Resnick stationary}, if the associated max-stable process $\eZ$ is stationary i.e., $
\eZ \fdd  \thhh  \eZ$  for any $h\in \TT$, where $L$ is the lag (backshift) operator and thus $ \thhh  \eZ(t)= \eZ(t- h),  h,t-h \in \TT.$\\
 For a positive $\sigma$-finite measure $\mu$ on $\mathcal{T}$, let
$\Pi_\mu=\sum_{i=1}^\IF \varepsilon_{(P_i,T_i)}$ be a PPP on $\R\times \TT$ with intensity
$e^{-p}dp \cdot  \mu(dt)$ being independent of anything else. If $Z$ is  a {\it Brown-Resnick stationary} and sample continuous Gaussian process on $\mathcal{T}= \R^d$, in view of \cite{DM}[Th.\ 2.1] \Ec{(see also \cite{DombryE}[Th.\ 2])} the following {\it Dieker-Mikosch representation}
\bqn{\label{dmrep}
\eZ(t) &\fdd&
\Ec{ \max_{i\ge 1}\Biggl( P_i  +Z_i(t-T_i) -\ln \int_{\mathcal{T}} e^{Z_i(s- T_i)} \, \mu(ds)   \Biggr)}, \quad
t\in \TT
}
is valid, provided that $\mu$ is a probability measure and a.s. $Z(\oldemptyset_\TT)=0$ with $\oldemptyset_\TT$ the origin of $\TT$. For notational simplicity hereafter we write simply $\oldemptyset$ instead of $\oldemptyset_\TT.$ \\
We shall show that \eqref{dmrep}
given in terms of the {\it tilted spectral processes} holds for general non-Gaussian  $Z$ and some positive $\sigma$-finite measure $\mu$ on $\TT$,   see \netheo{Th1} and \netheo{DMElef} below.
Motivated by \cite{DM} we present some useful conditions for the stationarity of $\eZ$.
\Ec{As a by-product we derive a new characterisation of Gaussian df's and give a new proof of Kabluchko's characterisation of Gaussian processes with stationary increments, see \netheo{korrK} and \netheo{korrP} in Section  2. }
 \item[B)]  An interesting class of stationary max-stable processes $\eZ(t),t\ge 0$ is constructed by Stoev in \cite{StoevSPA}, where $B(t),t\ge 0$ is a real-valued
 L\'evy process with $\E{e^{ B(1)}} < \IF $ and $Z$ is specified via \eqref{kanon}. We show
 in \netheo{kolmog} that for general $Z$ a two-sided extension of
$\eZ$ can be defined in terms of some {\it spectral process} $Y$ determined by $\Dhz, h\in \TT$.
\item[C)] If \eqref{zho} does not hold we modify the definition of $\Dhz$, see \nelem{lemDhV}.
Such a modification shows that the {\it tilted spectral processes} have  a component which is identifiable  and moreover  it determines the law of $\eZ$. Specifically, for any $t_1 \ldot t_n \in \TTT$ and $H$ the df of $(\eZ(t_1) \ldot \eZ(t_n))$, we derive the following (referred to as the {\it  inf-argmax representation})
\bqn{
\label{HR}
- \ln \Hxn &=& \sum_{h=1}^n e^{- x_h} \Psi_h(x) , \quad \forall x=(x_1 \ldot x_n) \inr^n,
}
with $\Psi_h$'s determined by the identifiable part of $\Dhz $,  see below \netheo{Euler}. \Ec{In the special case that $H$ is continuously differentiable
\eqref{HR} is a consequence of Euler's homogeneous functions theorem, see e.g., \cite{DombBook}[Eq. (9.11)]}.
\item[D)] The recent contribution \cite{SBK} introduced the generalised Pickands constant $\HZD$ defined by
\bqn{ \label{senti}
 \HZD= \limit{T} \frac 1 {T^{d}} \E*{ \sup_{t\in  \Zdd \cap [0, T]^d} e^{Z(t)}}, \quad \delta\ge 0,
 }
with the convention that $0 \Zd=\R^d$. 
We show in Section \ref{SecAppII} sufficient conditions that imply the positivity of  $\HZD$ for $\delta  \ge 0,d\ge 1 $ and derive further two new representations for  $d\ge 1$ and $\pk{Z(t)=- \IF}\ge 0, t\in \TT$ in terms of the so-called {\it spectral tail process} defined in \cite{BojanS}. Our findings for Pickands type constants are important since they are relevant for the calculation of the extremal index of multivariate stationary times series. 

\EIT

Organisation of the rest of the paper.  In Section 2 we establish the {\it Dieker-Mikosch representation} of $\eZ$
if  $Z$ satisfies \eqref{zh:finite} for any $h\in \TT$, and discuss further  some new conditions for the stationarity of $\eZ $. We continue with an application in Section 3 where we show how to construct a two-sided extension of $\eZ$. Section 4 is concerned with the general case that $Z$ takes values in $[- \IF, \IF)$. New formulas for $\HZD$   are  displayed   in Section \ref{SecAppII} followed by discussions and further results in Section \ref{SecDisc}. All the proofs are relegated to Section \ref{SecProof}.

\def\maxz{ \stackrel{\bigvee}{ \sim}}
\def\fI{ \mathcal{F}_{-\IF}}
\section{Max-stable processes with real-valued $Z$}
Let $Z, Z_i,i\ge 1,\Pi_\mu$ be as in the Introduction and suppose that for some  $h\in \TT$ the random variable (rv) $Z(h)$ satisfies  \eqref{zh:finite}  (the case $\pk{Z(h)=-\IF}> 0$ will be discussed in Section 4).  Let in the following 
$\mathcal{F}_{-\IF}$ denote the set of functions on $\mathcal{T}$ with values in $[-\IF,\IF)$ excluding the function $f$ equal to $-\IF$ and write $\mathcal{B}(\fI)$  for the $\sigma$-field generated by all evaluation maps. As in \eqref{defDhz} we define 
 \Ec{$\Dhz $ by $
	\Dhz(t)=  Z ^\HH (t)- \cE{Z ^\HH }(h),  t\in \TT$ with $Z ^\HH $ the exponential tilt of $Z$ by  $Z(h)$ i.e.,
$\pk{Z ^\HH  \in A}= \E{e^{Z(h)}\mathbb{I}\{Z \in A\}},\forall A \in \mathcal{B}(\fI).$} In view of \eqref{tDi}, if $\eta$ is a rv with values in  $[- \IF, \IF)$ being independent of $Z$ with  $\E{e^\eta}=1$, then
\bqn{\label{send}
\zeta_Y \fdd  \eZ,
}
with $\zeta_Y$ the max-stable process associated to $Y(t)= \eta+ Z(t), t\in \TT$.
 Although $Y$ and $Z$ are completely different  processes, we have that $\Xi_h Y \fdd \Dhz$.
 \Ec{Surprisingly, as shown below this fact holds for a general $Y$ satisfying  \eqref{send};  see also its extension in \nelem{prop2} covering the case $\pk{Z(h)= -\IF}>0$. }

\BEL \label{lem00} (\cite{HaanPickands}[Lemma 4.1,  A.A. Balkema]) If \eqref{send} holds  and a.s. $Y(h)=Z(h)=0$ for some $h\in \TTT$, then $ Y\fdd Z.$
\EEL

The claim of \nelem{lem00} is included  in  \cite{MolchanovBE} and \cite{StoevWangSPL};
 a direct proof is mentioned in \cite{FalkC} which is elaborated in \cite{FalkCH}[Lemma 1.1].
We present yet another proof in Section 7.

Since when $\pk{Z(h)>- \IF}=1$  we have $\Dhz(h)=0$ almost surely,  then \nelem{lem00} proves the uniqueness of $\Dhz$ (in the sense therein). This implies that $\Dhz$ can be determined directly in terms of $\eZ$.\\
Our next result  below confirms this. Moreover we show that  $\eZ$ possesses a {\it Dieker-Mikosch representation} determined by $\Dhz$ and \Ec{some  positive $\sigma$-finite measure $\mu$ on $\TT$,
 provided that $\SZ=\int_{\TT}  e^{Z(s)} \mu(ds)$ is a rv satisfying 
 		\bqn{\label{trift}
 			\pk{ \SZ< \IF}=1.
 		}
Note that the assumption $\E{e^{Z(t)}}=1,t\in \TT$ implies that  
\eqref{trift} holds for any  probability measure $\mu$  on $\TT$.} 
Throughout in the following  $H$ stands for the df of $(\eZ(t_1) \ldot \eZ(t_n))$ for some distinct $t_1 \ldot t_n\in \TTT$ and 
denote by $W=(W_1 \ldot W_n)$  an $n$-dimensional random vector with df $G$ given by
set 
\bqn{\label{dfG}
G( x)= \frac{1}{ \ln H(\oldemptyset  ) } \ln \Bigl( \frac{H(\min (x_1,0) \ldot \min(x_n,0)) }{H( x  )}\Bigr), \quad
 x=(x_1 \ldot x_n) \in [- \IF, \IF]^n .}
By the fact that $H$ is positively associated, see e.g., \cite{Molch,Res1987} we have $H(\oldemptyset)>0$. In view of \cite{Segers}[p. 278], see also \cite{SegersRoot}[Eq. (2.6)] the df's $G$ corresponding to different $t_i$'s are  the so-called generalised Pareto df's, here referred to as  the {\it associated GPD's} of $\eZ$. Set below 
$$(W_1^{(h)} \ldot W_n^{(h)}):=W^{(h)} \EQD  W \lvert (W_h>0 ), \quad h \in \{1 \ldot n\}$$ 
and note that $W_h^{(h)} $ is a unit exponential rv.

\def\oo{a}
\BT \label{Th1}  i) If $\pk{Z(h)> - \IF}=1$,  then  for distinct  $t_1=h, t_2 \ldot t_n \in \TTT$ 
\bqn{ \label{kro}
\Bigl(  \Dhz  (t_2) \ldot \Dhz  (t_n) \Bigr) \EQD  \Bigl(W_{2}^{(h)}- W_{1}^{(h)}\ldot W_{n}^{(h)} -W_{1}^{(h)}  \Bigr).
}
ii) If $\mu$ is a positive $\sigma$-finite  measure  on $\TT$ satisfying \eqref{trift} 
with $\Pi_\mu=\sum_{i=1}^\IF \varepsilon_{(P_i,T_i)}$ be a PPP on $\R\times \TT$ with intensity $e^{-p}dp \cdot  \mu(dt)$ being independent of anything else we have 
\bqn{ \label{RepNb}
\eZ(t) \fdd
\max_{i\ge 1} \Biggl( P_i  +  \Xi_{T_i} Z_i(t) -\ln \int_{\mathcal{T}} e^{  \Xi_{T_i} Z_i(s)   } \, \mu(ds) \Biggr),
 \quad t\in \TT.
}
\ET

\BEX \label{???xamp1}  Consider $Z(t)=B(t)- r(t,t)/2,t\in \TTT $ with $B$ a  centred, sample path continuous  Gaussian process with stationary increments  and covariance function $r$. Setting $\sigma^2(t)=r(t,t)$ we have
\BQN \label{dusig}
 \Xi_{R} Z(t)\fdd   B(t-R)- \sigma^2(t-R)/2, \quad t\in \TT
\EQN
for any real-valued rv $R$ independent of $B$. Hence
 \eqref{RepNb} reduces to \cite{DM}[Th.\ 2.1] 
 when $\mu$ is a probability measure.
\EEX

\BRM
If $\eZ(t),t\inr^d, d=1$ is stationary, and a.s. $Z(0)=0$, then by \eqref{dusig}  the  representation in \eqref{RepNb} agrees with the finding of \cite{SBK}[Th.\ 4].
\ERM

\def\TT{\mathcal{T}}
\def\vka{{0}}

Several contributions have investigated the stationarity of max-stable processes and particle systems, see for details \cite{HaanPickands, kab2009a, kab2009, kabluchko2011extremes, K2010, MR3188356, engelke14, MolchanovSPA, MolchKristin}. The main result of this section displays three criteria for the stationarity of $\eZ$. Below we define $L^b\Xi_a Y:=L^b(\Xi_a Y)$ by 
\bqn{ \label{abu} L^b \Xi_a Y(t)=
	 L^b(Y^{[a]}(t)- Y^{[a]}(a))=  Y^{[a]}(t-b)- Y^{[a]}(a)
}
for any $a,t-b,t\in \TT$ and some process $Y$ such that $\E{e^{Y(a)}}=1$ with $Y^{[a]}$ the tilted process by $Y(a)$. Recall that in our notation $\TT=\R^d$ or $\TT=\Z^d$ and $\oldemptyset$ is the origin in $\TT$.
\BT \label{ThM} Let $\eZ(t),t\in \TT$ be a max-stable process with unit Gumbel  marginals  and spectral process $Z$ defined via \eqref{W} and let for some $\sigma$-finite measure $\mu$ on $\TT$ 
the PPP $\Pi_\mu$ be as in the Introduction. If \eqref{zh:finite} holds  for any $h\in \TT$, then the following are equivalent: \\
 a) $\eZ$ is  stationary i.e., $\eZ \fdd \zeta_{ \thhh  Z}$ for any $h\in \TT$.\\
 b) For any positive $\sigma$-finite measure $\mu$ on $\TT$
  we have that \eqref{RepNb} holds with
$L^{T_i} \Xi_\oldemptyset Z_i$ instead of $ \Xi_{T_i} Z_i$, provided that $\SZ=\int_{\TT} e^{Z(t)} \mu (dt)$ is a positive finite rv. \\
 c) For any functional  $\Gamma:$  $\fI \to [0,\IF)$ which is
$\mathcal{B}(\fI) / \mathcal{B}(\R)$ measurable,  such that $\Gamma(f+c)= \Gamma(f),c\in \R,f\in [-\IF,\IF)^\TT $ holds, we have
 \bqn{ \label{thmc}
 \E[\big]{e^{Z(a+h)}\Gamma(Z)} = \E[\big]{e^{Z(a)}\Gamma( \thhh  Z)},
 \quad \forall h\in \TT,
 }
provided that the expectations exist. \\
d)  For any $a, a+h \in \TT $ 
 \bqn{\label{thmd}
  \Xi_{a+h} Z \fdd  L^h  \Xi_a  Z.
 }
\ET

\def\hh{ \thhh }

\def\inT{\in \TT}

\BRM  \label{remiii}
i) If $Z$ is as in \neexamp{???xamp1}, then statement $a) \Longrightarrow c)$ in \netheo{ThM} has  been shown in \cite{DM}[Lemma 5.2], whereas the non-Gaussian case is derived   in \cite{SBK}[Lemma 1] under the restriction that a.s. $Z(\oldemptyset)=0$. \\
ii) If statement c) and  d) in \netheo{ThM}  hold for any $h\in \TT$ and $a=\oldemptyset $ being the origin, then $\eZ$ is max-stable and stationary. \\
\ERM

\Ec{We present next two applications, a third one is displayed in Section 3.\\
	 Motivated by 	\cite{MR3454027}[Th.\ 1]  we derive below a new  characterisation of Gaussian distributions. Hereafter $(\cdot, \cdot)$ stands for the scalar product in $\R^d$. }\\
\Ec{\BT \label{korrK} Let $X$ be a $d$-dimensional random vector with non-degenerate components and define  $Z(t)= (t, X)  - \kappa(t),t \inr^d$, with $\kappa$ some measurable function satisfying $\kappa(\oldemptyset)=0$. Suppose that the associated max-stable process $\eZ(t), t\inr^d$ has unit Gumbel  marginals  and set  $\zeta_Z^\delta(t)=\eZ(t), t\in \delta \Z^d$. If for any $\delta \in (0,\IF), h\in \Zdd$
\bqn{ \Xi_h Z(t) \fdd  \thhh  Z(t), \quad t\in \Zdd,}
then $X$ is Gaussian $N_d( \mu, \Sigma)$ and $ \kappa(t)= (t, \mu) + ( t, \Sigma t)/2, t\inr^d$.
\ET
}

Our second application  is a different proof of Kabluchko's characterisation of Gaussian random fields with stationary increments 
stated in \cite{K2010}[Th.\ 1.1].

\BT \label{korrP} Let $B(t),t\in \R^d$ be a centred Gaussian process with non-zero variance function $\sigma^2$ such that $\sigma(\oldemptyset)=0$. The max-stable process $\eZ$ associated to $Z(t)=B(t)- \sigma^2(t)/2, t\in \R^d $ is stationary  if and only if $B$ has stationary increments.
\ET

\section{two-sided stationary max-stable processes}
Consider $\eZ(t), t\ge 0$ defined via \eqref{W}, where $Z(t)=  B(t)- t/2, t\ge 0$ with $B(t),t\inr$ a two-sided standard Brownian motion.
The seminal article \cite{bro1977} showed that $\eZ$ is max-stable and one-sided stationary. In view of \cite{kab2009}, in order to define $\eZ(t)$ also for $t< 0$  i.e., to define a two-sided stationary max-stable process $\eZ$, we can take $Z(t)= B(t)- \abs{t}/2, \forall t\inr$. This construction is fundamental since $B$ is both a centred Gaussian process with stationary increments and also a L\'evy process.
\COM{In the case that $B(t), t\ge 0$ is a general centred Gaussian process
with stationary increments, then $B(t), t\inr$ can be easily constructed since the covariance function of $B$ is defined by $\sigma^2(\abs{t})$ the variance of $B(t)$.   Moreover, by \cite{kab2009} the max-stable process
$\eZ(t),t\inr $ with $Z(t)= B(t)- \sigma^2(\abs{t})/2, t\inr$ is stationary.\\
}

Stoev showed in \cite{StoevSPA} that if $B(t),t\ge 0$ is a real-valued L\'evy process with Laplace exponent $\Phi(\theta)=
\E{e^{\theta B(1)}}$ being finite for $\theta=1$, then $\eZ(t), t\ge 0$ defined by \eqref{W} with $Z(t)= B(t)- \Phi(1) t$ is both max-stable and stationary. 
The recent contribution \cite{eng2014d} is primarily motivated by the question of how to define directly $Z(t), t< 0$ such that $\eZ(t),t\inr$ is both max-stable and stationary.  In Theorem 1.2 therein  a two-sided version of $Z$ and thus of $\eZ$ is constructed.  Specifically, as in \cite{eng2014d}  define $Z(t), t<0$ by setting
$Z(t)=  Z^{-}(-t), \quad t< 0,$ where  $Z^{-}(t), t\ge 0$ is independent of $Z(t),t\ge 0$ such that
$-Z^{-}(t)$ is the exponential tilt of $Z$ at $t$ i.e., in our notation since a.s. $Z(0)=0$,  then for any $t> 0$
$$ Z^{-}(t)= \Xi_{t} Z(-t)=  Z^{[t]}(0)- Z^{[t]}(t)= - Z^{[t]}(t).$$
Hence, in view of \cite{Kypri}[Theorem 3.9] (see also \cite{DeM15})
$Z^{-}$ is a L\'evy process with Laplace exponent  $\ln \mathbb{E}\{e^{\theta Z^{-}(1)}\}=\Phi(1- \theta)- (1- \theta) \Phi(1) .$
Our next result is not restricted to the particular cases of $Z$ being a L\'evy or a Gaussian process.

\BT \label{kolmog} Let $\eZ(t), t\ge 0$ be a max-stable and stationary  process determined by $Z$ as in \eqref{W} with $\E{e^{Z(t)}}=1, t\ge 0$. If  \eqref{zh:finite} holds for any $h \ge 0$,  then there exists a random process $Y(t),t\inr $ such that for distinct $t_1 \ldot t_n \inr $
\BQN\label{kolmogorov}
\Bigl ( Y(t_1)  \ldot Y(t_n)\Bigr) \EQD  \Bigl( \Xi_{h} Z(t_1+h) \ldot  \Xi_{h} Z(t_n+ h) \Bigr), \quad h:= - \min \Bigl(0, \min_{1 \le j \le n}t_j \Bigr)
 \EQN
and $\zeta_{Y}(t) \fdd \eZ(t),t\ge 0$. Moreover $Y(t),t\inr $ is  
{\it Brown-Resnick stationary}.
\ET
\COM{Note in passing that by \eqref{kolmogorov}, if  a.s. $Z(0)=0$, then for any $t\inr$
\BQN \label{dra}
 Y(t) &\EQD  &-Z^{[-t]}(-t)\mathbb{I}\{ t< 0\}+ Z(t)\mathbb{I}\{t\ge 0\}.
\EQN
}

\BEX \label{???xamp3} (Brown-Resnick process) Let $\eZ(t),t\ge 0$ be a max-stable process associated to $Z(t)=B(t)- \sigma^2(t)/2, t\ge 0 $ with $B(t),t\inr$
a centered Gaussian process with stationary increments and variance function $\sigma^2$. If $\sigma(0)=0$, by \neexamp{???xamp1}  and \eqref{kolmogorov} it follows easily that $Y \fdd  Z^* $, where $Z^*(t) = B^*(t)- \sigma^2(\abs{t})/2,t\inr$ with  $B^*(t),t\inr$
a centered Gaussian process with covariance function $(\sigma^2(\abs{t}) + \sigma^2(\abs{s}) - \sigma^2(\abs{t-s}))/2$.
Since $B^*$ has stationary increments, then by \netheo{korrP} $\eZ$ is stationary.
\EEX

\BEX \label{???xamp4}  (L\'evy-Brown-Resnick process) Suppose that $Z(t),t\ge 0$ is a L\'evy process with $\E{e^{Z(t)}}=1,t\ge 0$.
According to \cite{StoevSPA} the max-stable process $\eZ(t),t\ge 0$ associated to $Z$ is stationary.
Hence we are in the setup of \netheo{kolmog}, which ensures that $\zeta_{Y}(t), t\inr$ is a max-stable stationary extension of $\eZ$. Further for $s \le t < 0$ by \eqref{kolmogorov} and $Z(0 )=0$
\BQN \label{Wi}
 \bigl(Y(s)- Y(t), Y(t) \bigr) 
&\EQD  & \bigl( - Z^{[-s]}(t-s), Z^{[-s]}(t-s)- Z^{[-s]}(-s)\bigr).
\EQN
The assumption that  $Z(t),t\ge 0$ is a L\'evy process yields that $Y(s)- Y(t)$ is independent of $Y(t)$, and thus $Q(t)=Y(-t),t\ge 0$ has independent increments. If $s< 0 < t$, then also $Y(s)$ is independent of $Y(t)$. Further \eqref{Wi} implies that  $Q$ has  stationary increments.
\COM{Indeed, for $a,x$ positive
$$ W(x+a)- W(x)= Y(-x-a)- Y(-x) \EQD   - Z^{x+a}(a).$$
}
Since $Z$ has independent increments and $\E{e^{Z(x)}}=1$, then  for $a,x$ positive and $v\inr$
\BQNY
\pk{Q(x+a)- Q(x) \le v} 
= \pk{-Z^{[a]}(a)\le v}
 \EQNY
and thus $Y$ agrees with the definition of \cite{eng2014d}.
\EEX

\def\OO{ \Xi_h}

\section{General Spectral Processes}
In this section we assume that $Z(h)=- \IF$ for some $h\in \TT$ with non-zero probability.
Write next  (set below $0 \cdot \IF=0$)
\bqn{ \label{zH}
Z \fdd J_h  V_h + (1- J_h) W_h, 
}
where $J_h$ is a Bernoulli rv with
$$\pk{ J_h =1}= \pk{Z(h)> - \IF} \in (0,1]$$
and $$V_h \fdd Z \lvert (Z(h)> -\IF), \quad  W_h \fdd Z \lvert (Z(h)=- \IF).$$ Furthermore, $J_h, V_h, W_h$ are mutually independent and
$$ \pk{V_h(h)> - \IF}= \pk{W_h(h)= - \IF}= 1.$$
For  $V_h ^\HH $ given via exponential tilting as
$$\pk{V ^\HH _h \in A}= \E{ e^{ V_h(h) - \ln \E{e^{V_h(h)}}}  \mathbb{I} \{ V_h \in A\}}, \quad  A\in  \mathcal{B}(\fI)$$
define the {\it tilted spectral process} $\Dhz $ by
\BQN \label{spectWroc}
\OO Z (t)=   J_h  \Theta_h (t) +  (1- J_h) [W_h(t)- V ^\HH _h(h)]- \ln \pk{J_h=1},  \quad t \in \mathcal{T},
 \EQN
where 
\bqn{\label{ThetH} \Theta_h(t):= \Xi_h V_h(t) = V ^\HH _h(t)- V ^\HH _h(h).
}	
 We shall consider $V^\HH_h$ to be independent of $J_h$ and $W_h$. \\
  The next result establishes the counterpart of \eqref{repA}. Further, we give a representation of $\eZ$ which is motivated by \cite{KabExt}[Th.\ 2].

\BEL  \label{lemDhV} For any $h\in \TT$ we have
$\eZ \fdd \zeta_{\Dhz }$. Moreover, for any probability measure $\mu$ on $\TT$ we have
$\eZ \fdd \eta$ where $\eta(t)= \max_{i\ge 1} (P_i+ \Xi_{T_i} Z_i(t)), t\in \TT$ with $(P_i, T_i)$'s the points of a PPP on $\R \times \TT$ with intensity
$e^{-p} dp  \cdot \mu(dt)$ being independent of $Z_i,i\ge 1$.
\EEL

In view of \nelem{prop2} in Appendix  $ \Theta_h,h\in \TT$ is the identifiable part of the family  of {\it tilted spectral processes} $\Dhz, h\in \TT$. Moreover, as shown below $ \Theta_h, h\in \TT $ determines the law of $\eZ$. 

\BT \label{Euler} (\underline{Inf-argmax representation}) For any distinct $t_i \in \TT,  i\le n$  the
df $H$ of    $(\zeta_Z(t_1) \ldot \zeta_{Z}(t_n))$ is given by
\BQN\label{Euler:rep}
- \ln  \Hxn &=&\sum_{k=1}^n e^{- x_k} \Psi_k(x ), \quad \text{with }
\Psi_k(x)= \pk[\Big]{ \inf {\rm argmax}_{1 \le i \le n} \big (\Theta_{t_k}  (t_i) - x_i\big )=k}
 \notag
\EQN
for any $x=(x_1 \ldot x_n)\inr^n$.
\ET
We conclude this section with an extension of \netheo{ThM}.
 
\BT \label{ThMB} Let  $Z(t),t \in \TT$ be a random process with values in $[- \IF, \IF)$.
If $\E{ e^{Z(t)}}= 1,  t \in \TT $ and $\eZ$ is given by \eqref{W}, then
the following are equivalent: \\
 a) $Z$ is {\it Brown-Resnick stationary}.\\ 
 b) For any $\Gamma$ as in the statement c) in \netheo{ThM} 
  \bqn{\label{xYx}
\E*{ e^{Z(a+h)} \Gamma(Z)}& = 
\E*{ e^{Z(a)} \Gamma(L^h Z)}=
\E[\big]{ \Gamma( L^h\Theta_{a}  )}, \quad a,a+h\in \TT.
  	 }
c)  For any $a, h\in \TT$ and $Z$ with representation \eqref{zH} we have
 \bqn{\label{thmdBB}
  \Theta_{a+h} \fdd L^h \Theta_a.
 }

\ET
Otherwise specified,  hereafter we set 
$$ \Theta= \Theta_\oldemptyset, \quad  x=(x_1 \ldot x_n).$$
{\BRM
If $\Gamma, \eZ$ are as in \netheo{ThMB}, then \eqref{xYx} is equivalent with 
\bqn{ \label{pla}
	\E*{ \mathbb{I}\{Z(-h)>- \IF\}e^{Z(-h)} \Gamma(L^{h}Z)}& = 	\E[\big]{ \Gamma( \Theta  )}, \quad h\in \TT.
}
Hence the {\it inf-argmax representation} in \eqref{Euler:rep} simplifies to 
{\BQN\label{Euler:rep3}
	- \ln  \Hxn	&=&\sum_{k=1}^n e^{- x_k} 
	\pk[\Big]{ \inf {\rm argmax}_{1 \le i \le n} \big (L^{t_k} \Theta  (t_i) - x_i\big )=k}, 
	\quad x\inr^n
	\EQN}
and thus we conclude that the fidi's of $\eZ$ are given in terms of those of  $\Theta$.
\ERM

 	\section{Generalised Pickands Constants} \label{SecAppII}
 	Given $Z(t),t\in \R^d, d\ge 1$ with representation \eqref{kanon} we  define  for any $\delta >0$
 	the generalised Pickands constant $\HZD$ as in \eqref{senti} i.e.,
\bqny{
	\HZD= \limit{T} \frac 1 {T^{d}} \E[\Big]{ \sup_{t\in  \Zdd \cap [0, T]^d} e^{Z(t)}}.
}
 	A canonical example here is the Brown-Resnick stationary case with $Z(t)= \sqrt{2} B_\alpha(t) - \abs{t}^\alpha, t\inr,$ where  $B_\alpha, \alpha \in (0,2]$ is a standard fractional Brownian motion with Hurst index $\alpha/2 \in (0,1]$. For this case  $\HZD$ is the classical Pickands constant, see e.g., \cite{PickandsB, Pit96, demiro2003simulation, Krzys2006Pickands,mi:17, DeHaJi2016} for its properties.\\
 	The recent contribution \cite{SBK} investigates $ \HZD$ under the assumption that a.s. $Z(\oldemptyset)=0$  and $d=1$. 
 	In this section we shall assume that that  $\eZ(t),t\in \TT$ 
 	 is max-stable, stationary and has unit Gumbel   marginals. In order to show the positivity of $ \HZD$,  we shall assume further that
 	\bqn{\label{mixing}
 		\pk*{ \int_{ \R^d} e^{Z(t)} \lambda(dt)< \IF}=1,
 	}
 	where $\lambda$ is the Lebesgue measure on $\R^d$. In  light of \cite{dom2016}[Th.\ 2] (see also \cite{WangStoev}) \eqref{mixing} is
 	equivalent with
 	\bqn{ \label{limNorm} 
 		\pk[\Big]{ \limit{\norm{t} }  Z(t)=- \IF}=1.
 		}
 	 Under \eqref{mixing}, as in \cite{SBK}[Th.\ 1], if a.s. $Z(\oldemptyset)=0$, then for any $\delta>0$ \Ec{and $d=1$}
 	   \BQN \label{JohannaSKBE}
 	\HZD =  \mathbb E\left\{\frac{\sup_{t\in \Zdd  } e^{ Z(t)}}
 	{ \delta^d \int_{  \R^d  } e^{Z(s)}\mu_\delta(ds) } \right\} \in (0,\IF),
 	\EQN
where $\mu_\delta$ denotes the counting measure on $\Zdd$. 
 If $\pk{Z(\oldemptyset) =0}<1$  the expression in \eqref{JohannaSKBE}  needs to be modified, since by \nelem{lemD} in Appendix, for any $T>0$ and $d\ge 1$ we have 
 	\bqn{\label{newiPick}
 		\HZD &=&  \limit{T}
 		\int_{ [0,1]^d} \mathbb E\left\{\frac{\sup_{t\in   \Zdd \cap [-h T,(1-h)T]^d } e^{  \Theta (t)}}
 		{ \int_{ \Zdd \cap [-hT, (1-h)T]^d} e^{ \Theta(s)}\mu_\delta(ds) } \right\}\mu^T(dh) \notag\\
 		& =:&  \limit{T}
 		\int_{ [0,1]^d}  \eta_T(h) \mu^T(dh), \quad 
 	}
where $\mu^T (dh)= \mu_\delta ( T dh)/T^d.$ In applications, often Pickands-type constants corresponding to $\delta=0$ appear. In order to define $\mathcal{H}_Z^0$, we shall  suppose further that $\eZ$ has cadlag sample paths. This is equivalent with $\E{ \sup_{t\in K} e^{Z(t)}}< \IF$ for any compact set $K \subset \R^d$, see \cite{kabDombry}.
The definition of $\mathcal{H}_Z^0$ is exactly as in \eqref{6.6} where we interpret $0 \Z^d$ as $\R^d$ i.e., $\mathcal{H}_Z^0 = \limit{T} \frac 1 {T^{d}} \E{ \sup_{t\in [0, T]^d} e^{Z(t)}}$. The existence and the finiteness of $\mathcal{H}_Z^0$ follow easily by the stationarity of $\eZ$. As in the case $\delta>0$ of interest is the positivity of
$\mathcal{H}_Z^0$, alternative formulas or bounds. For $Z$ Gaussian or L\'evy and $d=1, Z(0)=0$ a.s. \cite{SBK}  shows that under some weak restrictions
$$ \mathcal{H}_Z^0 = \mathbb E\left\{\frac{\sup_{t\in \R^d  } e^{ Z(t)}}
{  \int_{  \R^d  } e^{Z(s)}\lambda(ds) } \right\}\in (0,\IF).$$
Since the aforementioned results in the literatute cover only the case $Z(0)=0$ a.s.\ and $d=1$, below we shall derive two formulas for $\HZD, \delta>0$ for the general case $d\ge 1$ and $\pk{Z(\oldemptyset)=- \IF} \ge 0$. If $\delta=0$, we give  a positive lower bound for $\mathcal{H}_Z^0$. 
\BT \label{mentalist} Let  $Z(t),t\in \R^d,d\ge 1$ be such that the associated max-stable process $\eZ$ is stationary with unit Gumbel marginals and  suppose that \eqref{mixing} holds.\\
 	i) For any $\delta >0$
 	\bqn{\label{newiPick2}
 		\HZD
 		&=& 		\mathbb E\left\{\frac{\sup_{t\in \Zdd  } e^{  \Theta (t)}}
 		{ \delta^d \int_{  \R^d  } e^{ \Theta(s)}\mu_\delta(ds) } \right\} =:C^\delta \in (0,\IF).
 	}
 	Moreover, if further the fidi's of $ \Theta$ are absolutely continuous, then 
 	\bqn{\label{neformP}
 		\HZD = \frac{1}{\delta^d} \pk[\Big]{ \sup_{t\in \Zdd}  \Theta(t)=0}.
 	}
 ii) We have for $\delta=0$ 
  	\bqn{\label{newiPick3}
 	\mathcal{H}_Z^0
 	&\ge & 		\mathbb E\left\{\frac{\sup_{t\in \R^d  } e^{  \Theta (t)}}
 	{ \int_{  \R^d  } e^{ \Theta(s)}\lambda(ds) } \right\} =:C^0\in (0,\IF).
 }	
\ET
  	
 	\BRM  \label{D_remDelta0}
i) If  $Z$ is as in \neexamp{???xamp3}, then $Z=  \Theta$ and thus  \eqref{neformP} follows from \cite{SBK}[Th.\ 1] combined with  Theorem 8 and Remark 9 in \cite{KabExt}. A direct proof for $B$ being a standard fractional Brownian motion is given in \cite{DiekerY}[Prop. 4]. The lower bound $ 	\mathcal{H}_Z^0 \ge C^0$ is derived in \cite{SBK}[Th.\ 1] for $d=1, Z(0)=0$. It is of interest (and open question) to know general tractable conditions that yield $ 	\mathcal{H}_Z^0=C^0$. \\
ii) If $Z(t),t\inr^d$ is as in \netheo{mentalist}, then for any $\delta\ge 0$ we have
 the following Gumbel limit theorem (see for related results \cite{DM,stoev2010max,MIKYuw})
 \bqn{\label{GL}
 	\limit{T} \pk[\Big]{\sup_{t\in \Zdd \cap [0,T]^d} \eZ(t) \le x+ d\ln T}= e^{-\mathcal{H}_Z^\delta e^{-x} }, \quad \forall x\inr.
 }
 \Ec{Consequently, $\delta^d\mathcal{H}_Z^\delta>0,\delta>0 $ is the extremal index of the stationary random field $\eZ(t), t\in \Zdd $.} \\
 iii) \Ec{Generalised Pickands constants have appeared also in the non-Gaussian setup, see e.g., \cite{KW}[Lemma 5.16]. Recently,
 	\cite{WangOrstein} presented a Pickands type constant arising in the connection with semi-min-stable processes.}
 \ERM

 \COM{
 	Suppose that $Z$ have  continuous sample paths (similar arguments hold if the sample paths are assumed to be cadlag).
 	Let thus $E=C(\R^d)$ be the space of continuous real-valued functions $f:\R^d \to \R$.  The assumption that
 	the max-stable process  $\eZ$ has continuous sample paths, which is equivalent  (see \cite{kabDombry}) with $\E{\sup_{t\in K} e^{Z(t)}} < \IF $ for any compact $K\subset \R^d$. In this section we shall assume that the latter condition holds.

 	Stationary max-stable processes $\eZ$ with continuous sample paths  are important for the asymptotic theory of Gaussian random fields. Indeed, since we assume that $\eZ$ has continuous trajectories for any $T>0$ we have
 	\bqn{\label{eS}
 		\pk{\sup_{t\in [0,T]^d} \eZ(t) \le x+ d\ln T}= e^{-  e^{-x}  T^{-d} \mathcal{H}_Z[T]}, \quad \forall x\inr,
 	}
 	where $\mathcal{H}_Z[T]= \E{\sup_{t \in [0,T]^d} e^{Z(t)}}$. By the stationary of $\eZ$, we have that $\mathcal{H}_Z[T],T\ge 0$ is subadditive, and thus the following Gumbel limit theorem (see for related results \cite{DM,stoev2010max,MIKYuw})
 	\bqn{\label{GL}
 		\limit{T} \pk{\sup_{t\in [0,T]^d} \eZ(t) \le x+ d\ln T}= e^{-x \mathcal{H}_Z}, \quad \forall x\inr,
 	}
 	with $\mathcal{H}_Z= \limit{T} T^{-d} \E{\sup_{t \in [0,T]^d} e^{Z(t)}}\in [0,\IF)$. The study of stationary $\eZ$ is important for the asymptotic theory of Gaussian random fields, since if $Z$ is a {\it Brown-Resnick stationary} Gaussian random field, then $\mathcal{H}_Z$ is the generalised Pickands constant, see \cite{DM,DiekerY, SBK,MIKYuw, stoev2010max,debicki2002ruin}.\\
 	
 	It is well-known that  $C(\R^d)$ can be equipped with a metric which turns it into a  Polish space (use the metric of the uniform convergence on compact sets).
 	Furthermore,  the $\sigma$-field $\mathcal{B}(E)$ generated by the open balls of that metric agrees with the Borel $\sigma$-field. Hence the equality in law of two random fields with paths in $C(\R^d)$ is equivalent with the equality in law of the random fields with respect to the Borel $\sigma$-field $\mathcal{C}$. Consequently, if $\eZ$ has continuous sample paths, we define $Z ^\HH  $ to have also continuous sample paths by writing its law as $\pk{Z ^\HH  \in A}= \E{e^{Z(h)} \mathbb{I}\{Z\in A\}}$ for any $A \in \mathcal{C}$.
 	We have by the above findings that
 	\BQN  \pk{ \eZ \in A} = \pk{ \zeta_{\Dhz} \in A}, \quad \forall h\in \R^d, A \in \mathcal{C}.
 	\EQN
 	We state next the following reformulation of \netheo{ThM} for this special case.
 	
 	\BT \label{ThM1} Let the random field $Z(t),t\inr^d$ be as in \netheo{ThM}, and suppose that $Z$ and $\eZ$ have
 	continuous sample paths. The following are equivalent: \\
 	a) $\eZ$ is a stationary max-stable random field and the Gumbel limit theorem in \eqref{GL} holds. \\
 	b) For any arbitrary probability measure $\mu$ on $\R^d$ 
 	$$ \pk{ \eZ \in A} = \pk{ \zeta_{T,Z} \in A}, \quad \forall A \in \mathcal{C},$$
 	where $\zeta_{T,Z}$ is as in \netheo{ThM}.\\
 	c) For any positive functional $\Gamma:$  $C(\R ^d) \to \R$ which is
 	Borel measurable (i.e., $\mathcal{C} / \mathcal{B}(R)$ measurable), such that $\Gamma(f+c)= \Gamma(f),f\in E$ holds for any constant function $c$, then
 	\bqn{ \E{e^{Z(a+h)}\Gamma(Z)} = \E{e^{Z(a)}\Gamma( \thhh  Z)}= \E{\Gamma( \Xi_{\oo}( \thhh Z))}, \quad \forall a,h\inr^d,
 	}
 	provided that the expectations exist. \\
 	d)  For any $a,h\inr^d$
 	\bqn{\label{thmd}
 		\pk{  \Xi_{a+h} Z \in A} = \pk { \Xi_a ( \thhh  Z) \in A}, \quad \forall A \in \mathcal{C}.
 	}
 	\ET
 	
 }

 \section{Discussions \& Further Results}  \label{SecDisc}
 \def\N{\mathbb{N}}
 \def\hh{ ^\HH \!}
 \def\TT{\mathcal{T}}

 	\subsection{Tilted processes} If $Z(t),t\in \TT$  is a random process with finite 
 $\pF= \ln \E{e^{Z(h)}},h\in \TT$, then we define $Z ^\HH $ by
 	\bqny{ \label{Fubini}
 		\pk{Z ^\HH \in A}&=& \E{ e^{Z(h)- \pF}   \mathbb{I}\{Z \in A\}}, \quad
 		A \in \mathcal{B}(\fI).
 	}
The exponential tilting of df's in the exponential family can be calculated explicitly. In particular, for the Gaussian case, the tilted process is again Gaussian, with the same covariance function, but modified mean,  see \cite{glyn}[p. 130] or \Ec{\cite{DombryE}[Lemma 1]}.

 	\BEL \label{tiltGauss}
 	Let $Z(t),t\in \mathcal{K}$ be a Gaussian process defined on some arbitrary parameter set $\mathcal{K}$
 	with covariance function $r$. For any $h\in \mathcal{K}$ the random process $Z ^\HH $ 
 is again Gaussian and moreover
 	\bqn{ Z ^\HH (t)\fdd Z(t)+ r(h,t), \quad t \in \mathcal{K}.
 	}
 Conversely, if $Z^\HH(t), t \in \mathcal{K}$ is for some $h\in \mathcal{K}$ 
 a Gaussian process with covariance function $r$ and mean $r(h,t)$, then $Z$ is a centred Gaussian process with covariance function $r$.
 	\EEL

 \BEX \label{???xamp1B}  Consider $Z(t)=B(t)- r(t,t)/2,t\in \mathcal{K} $ with $B$ a  centred Gaussian process with covariance function $r$. For any $h \in \mathcal{K}$ by \nelem{tiltGauss} $Z ^\HH (t)\fdd  B(t)- r(t,t)/2 + r(t,h)$ implying that
 \bqn{\label{exDhZ}
 	\Dhz (t) = Z ^\HH (t)- Z ^\HH (h)\fdd   B(t)- B(h) - Var(B(t)- B(h))/2, \quad t\in \mathcal{K},
 }
 which agrees with the definition of $\Dhz $ given in \eqref{defDhz}.
 \EEX

 \BEX \label{EulerHR}  If $\eZ$ is as in Example \ref{???xamp3}, then by \eqref{exDhZ}
 $$ \Theta_{t_k}(t)=  \Xi_{t_k} Z(t) \fdd B(t)- B(t_k)- Var(B(t)-B(t_k))/2 
 , \quad t,t_k\in \TTT,
 $$
 then  \eqref{HR}  holds with  $\Psi_h$ the df of the  Gaussian random vector $(\Theta_k(t_1) \ldot \Theta_k(t_n))_{-h}$ (the subscript $_{-h}$ means that the $h$th component is dropped).  Such a  representation of max-stable H\"usler-Reiss df has been derived by another approach in \cite{Joe92}, see also \cite{Joe}.

 \EEX

In order to calculate $\Xi_h Z$ when $Z(h)>- \IF$ almost surely, one can use alternatively \eqref{kro}. If $\pk{Z(h)> - \IF}> 0$, then in view of \nelem{prop2}, we have that \eqref{kro} holds with $\Theta_h$ insetead of $\Xi_h Z$. Hence when the distributions of the associated GPD's of $\eZ$ are known, we can calculate $\Theta_h$ using the right-hand side of \eqref{kro}.

 	\def\OD{{\Xi}}
\def\ODH{\OD_h}

 	\subsection{Dieker-Mikosch representation}
In view of our findings in Section 2 we have the representation (recall \eqref{RepNb})
\bqn{\label{reifw} \eZ(t) \fdd \max_{i\ge 1}\bigl (P_i+ F_{i}(t,T_i) \bigr), \quad t\in \TT,}
where
$(P_i,T_i)$'s are the points of a PPP $\Pi_\mu$ on $\R \times \TT$ with intensity $e^{-p} dp \cdot \mu(dt)$ being further independent of $F_{i}$'s which are independent copies of a random shape function $F$ defined by
\bqn{\label{shapeF}
F(t,h)= \Xi_h Z(t) -\ln \int_{\mathcal{T}} e^{  \Xi_h Z (s)   } \mu(ds), \quad h,t\in \TT,
}
with  $\mu$  a positive $\sigma$-finite measure on $\TT$ (recall $Z_i$'s are independent of the points of $\Pi_\mu$).\\
Next, we shall assume that $\mathcal{M}_Z= \sup_{t\in \TT} e^{Z(t)}$ and $ \mathcal{S}_Z = \int_{\TT} e^{Z(t)} \mu(dt) $  are non-negative and finite rv's.  The representation \eqref{RepNb} of $\zeta_Z$ is shown under the assumption that
 	$\mathcal{S}_Z  $ is a.s.\ positive. Since we assume that $\E{e^{Z(t)}}=1, t\in \TT$,
 	if $\mu$ is a probability measure,  then  the finiteness of $\mathcal{S}_Z $ is guaranteed also for general {\it spectral processes} $Z$ with values in $[- \IF, \IF) $.  However, $\mathcal{S}_Z $ can be equal to zero with non-zero probability. Therefore, in this section the {\it Dieker-Mikosch representation} for $Z$ with values in $[- \IF, \IF) $  will be shown under the following restriction
 	\bqn{\label{armpush}
 		\pk{\mathcal{M}_Z>0,\mathcal{S}_Z  = 0}=0.
 	}
If $\mu$ possesses a positive probability mass function $p(t), t\in \TT=\{ t_1 \ldot t_n\}$, then
\bqny{
	\pk{\MZ>0,\SZ  = 0}
	= \pk[\Big]{\MZ>0,  \sum_{k=1}^n p(t_k )e^{  Z(t_k)}=0}=0 ,
}
hence \eqref{armpush} is valid for such $\mu$. Similarly, \eqref{armpush} holds also for $\mu$ the counting measure on $\TT= \Z^d$ and we do not need further conditions to show that $\mathcal{S}_Z$ is a rv.
 	\BT \label{DMElef}
 	Let $\eZ(t), t\in \TT$ be a max-stable process with representation \eqref{W} and measurable {\it spectral process} $Z$ satisfying \eqref{zH}. If $\mu$ is
 a positive $\sigma$-finite measure on $\TT$ such that \eqref{armpush} is valid with  $\MZ,\SZ$ being two non-negative finite rv's, then \eqref{reifw}
holds with random shape function $F$ given by
\bqn{ \label{DMB}
	F(t,h)= \ThetaH(t) -\ln \int_{\mathcal{T}} e^{  \ThetaH(s)   } \mu(ds), \quad h,t\in \TT.
}
\ET

 	{
 		\BRM  
If $\mu$ is a probability measure on $\TT$ and $T$ independent of $\Theta$ has law $\mu$, then under the assumptions of
 		\netheo{DMElef}
 		\bqn{\label{frep}
 			- \ln  \Hxn  
 				&=& \E*{  \max_{1 \le k \le n} e^{- x_k+ \Theta_T  (t_k)  - \ln \bigl(
 						\int_{\TT} e^{  \Theta_T (s)} \mu(ds) \bigr) }},
 					\quad x\inr^n.
 			}
 		Based on  \eqref{frep} the  simulation method developed in \cite{DM} can be applied to the general case of $Z$ as shown in \cite{DombryE}.
		\ERM
{ 
	For the case that $Z$ is  {\it Brown-Resnick stationary} we have the following mixed-moving-maxima (M3) Dieker-Mikosch representation:

\BK\label{9:20ora} Under the assumptions of \netheo{DMElef}, if further $\eZ$ is stationary, then \eqref{reifw} holds with random shape function $F$ given by 
\bqn{ \label{DMB2}
	F(t,h)= L^h \Theta (t) -\ln \int_{\mathcal{T}} e^{ L^h \Theta(s)   } \mu(ds), \quad h,t\in \TT.
}
Conversely, if for some  random process $\Theta(t) ,t\in \TT$ 
we have that $ \int_{\mathcal{T}} e^{ \Theta(s)   } \mu(ds)  $ is a positive rv
with $\mu$ the Lebesgue measure on $\TT$, then $\zeta(t),t\in \TT$ with representation given by the right-hand side of \eqref{reifw} and random shape function $F$ given in \eqref{DMB2} is max-stable and stationary. 
\EK
}

\BRM  In the special case  $d=1$ and $\pk{Z(0)=0}=1$  the representation  \eqref{DMB2} is 
stated in \cite{SBK}[Thm. 3.1].
\ERM

 \subsection{Max-stable processes with \FRE  marginals}
In various applications max-stable processes $\eZ(t),t\in \TT$ with \FRE   marginals are considered, see e.g., \cite{PadoanA, PadoanB,Robert0,Opitz, MR3493175,raey}. Specifically, we define
\bqn{\label{Pf}
	\eZ(t) \fdd \max_{i\ge 1} P_i Z_i(t), \quad t \in \TTT,}
where $\sum_{i=1}^\IF \delta_{P_i}$ is a  PPP on $(0,\IF)$ with intensity $x^{-2} dx$ being independent of $Z_i,i\ge 1$ which are independent copies of  a non-negative random process $Z(t),t \in \TT $ with $\E{Z(t)}=1, t\in  \TT $.
Let $V_h \fdd Z \lvert (Z(h)> 0)$; recall that $1=\E{Z(h)}= \E{V_h(h)} \pk{Z(h)>0}$. As in the Gumbel case, the {\it tilted spectral processes} 
$\ThetaH:= \Xi_h V_h, h\in \TT$ are  defined by (interpret below $0/0$ as $0$)
\bqn{\label{6.6} 
	\pk{\ThetaH \in A}&=&
	\E[\Big]{ \frac{ V_h(h)}{ \E{ V_h(h)}} 	\mathbb{I}\{ V_h/V_h(h)\in A\}} 
= 	 \E[\Big]{ Z(h)\mathbb{I}\{Z(h) > 0\}		\mathbb{I}\{ Z/Z(h) \in A\}} \notag \\
	&=& \int_{\TT} f(h)  \mathbb{I}\{f(h) > 0\}
	\mathbb{I}\{ f/f(h) \in A\}   \nu( df)
	, \quad   A \in \mathcal{B}(\mathcal{F}_0),
}
where  $\mathcal{F}_0$ is the set of non-negative functions on $\TT$ excluding the zero function endowed with $\sigma$-field $\mathcal{B}(\mathcal{F}_0)$ and $\nu$ stands for  the law of $Z$. \\
If  $H$ denotes the df of $(\eZ(t_1) \ldot \eZ(t_n))$, then  its marginals are unit \FRE and  moreover  its {\it inf-argmax representation} is given by }
\bqn{ \label{HR2B}- \ln \Hxn 	&=&
	\sum_{k=1}^n \frac 1 {x_k} \pk[\Big]{ \max_{1 \le i < k}\Theta_{t_k}(t_i) < \frac{ x_i}{x_k},
		\max_{k   \le i < n }\Theta_{t_k}(t_i)   \le  \frac{ x_i}{x_k}}\notag\\
	&=& 	\sum_{k=1}^n \frac 1 {x_k} \pk[\Big]{\inf {\rm argmax}_{1 \le i \le  n} \Bigl( \frac{\Theta_{t_k}(t_i)  }{x_i}\Bigr) = k}, \quad x\in (0,\IF)^n.
}
Note in passing that if $\eZ$ is stationary, then by \eqref{thmdBB}  $$\Theta_{t_k} \fdd L^{t_k} \Theta ,
$$
hence \eqref{HR2B} is determined only by $\Theta$ in this case 
(recall that we set $\Theta=\Theta_\oldemptyset$). \\
 Next, as in \cite{dom2016} define the {\it functional PPP} $\Phi=\{ \phi_i, i\ge 1\}$ on $\mathcal{F}_0$ with $\phi_i= P_i Z_i, i\ge 1$ and intensity measure $q$ determined by 
\bqn{ \label{mua} 
	q(A)= \int_{A} q(df)=\int_0^\IF  \pk{u Z \in A} 	u^{-2} du
	, \quad A \in \mathcal{B}(\mathcal{F}_0).
}
{
We have that $H$ is determined by $q$ as follows 
	\bqn{- \ln H(x)= q\{f\in \mathcal{F}_0: f(t_i)> x_i, \text{ for some } i=1 \ldot  n   \}, \quad x \in (0,\IF)^n.
		\label{avi}
	}
	\begin{remark} \label{remQ}
		A direct implication of \eqref{mua},\eqref{avi} and \nelem{lem00} is that for any two random processes $Y(t), Z(t),t\in  \TT$ such  that $\E{Z(t)}\in (0,\IF),t\in \TT$ and  
\bqn{ \label{ajv} 
	\int_0^\IF \pk{ u Z \in A}u^{-2} du = \int_0^\IF \pk{ u Y \in A}u^{-2} du, \quad A \in \mathcal{B}(\mathcal{F}_0), 
}
		then  $Z \fdd Y,$ provided that almost surely $Z(h)=Y(h)=1$ for some $h\in \TT$. Note further that in view of \nelem{prop2} the relation \eqref{ajv} implies $ \ThetaH \fdd 
		\Xi_h U_h$ with $U_h \fdd Y \lvert( Y(h)> 0)$. 	
	\end{remark}	
	Denote by $\phi^+_h,h\in \TT$ the extremal function at $h$ (defined in \cite{MR3024101,DombryE}) i.e., this is the set of functions $\phi \in \Phi$ such that $\phi(h)=Z(h)$. 	The next result includes two key findings of \cite{DombryE}[Prop. 1 \& 2].  
	
	\BT\label{minior} i) $\Phi_h=\Phi \cap \{ f \in \mathcal{F}_0: f(h)>0\}, h\in \TT$ is a  PPP  with intensity 
	\bqn{q_h(A)=\int_{A} \mathbb{I}\{ f(h)>0 \}q(df) = \int_{0}^\IF \pk{ u \ThetaH  \in A}u^{-2} du,  
		\quad A\in \mathcal{B}(\mathcal{F}_0) . \label{eq:min} 
	}
	ii) $\ThetaH ,h \in \TT$ is equal in law with  $\phi^+_h/Z(h)$ which is independent of $Z(h).$\\
	iii) Let $\mathcal{V}(t),t\in \TT$ be a random process such that $\mathcal{V}(h)=1$ for some $h\in \TT$. If  
	\bqn{
		q_h(A)
		= \int_0^\IF  \pk{ u \mathcal{ V}  \in A}u^{-2} du,    \quad A\in \mathcal{B}(\mathcal{F}_0),
	}	
	then $\mathcal{V}$ has the same law as $\ThetaH$.
	\ET
		\BRM 
In view of \netheo{minior} $ii)$ the stationarity of $\eZ$ is equivalent with 
	$\phi^+_h/Z(h),h\in \TT $ has the same distribution as $L^h \phi^+_0/Z(0), h\in \TT$  	(here $L^h \phi^+_0$ is the shifted extremal function at 0 corresponding to $ \{ P_i L^h Z_i, i\ge 1 \} $). 	Consequently,   \netheo{ThMB} statement $d)$ follows by the stationarity of $\eZ$\footnote{This short proof was kindly suggested by a reviewer.}. \\
	Alternatively, by \netheo{minior} $i)$ we have that the stationarity of $\eZ$ implies that $\Phi_h$ has the same law as $L^h \Phi_0= 
	\{ P_i L^h Z_i, i\ge 1 \} \cap \{ f\in \mathcal{F}_0: f(0)>0\}$  for any $h\in \TT$,  which  by $i)$ and $iii)$ in the aforementioned theorem yields that statement $a)$ implies  	statement $d)$ in \netheo{ThMB}.\\
\ERM
\subsection{Tilt-Shift Formula}
Let  $X(t),t\in \Z^d$ be a real-valued stationary time series. Commonly, $X$  is called 
jointly regularly varying  with index $\alpha>0$, if  the random vectors $(X(t_1) \ldot X(t_n)), t_i \in \Z^d, i\le n$ are for any $n\in \N$ regularly varying with index $\alpha$. For such $X$, as shown in \cite{BojanS} there exists the so-called {\it spectral tail process} (STP)
$\Theta(t),t\in \Z^d$ with $\Theta(0)=1$ a.s.\  that satisfied the {\it time-change formula}, which in our language translates to 
\bqn{\label{tcf}
	\E*{ \mathbb{I}\{\abs{\Theta(-h)}>0\} \abs{\Theta(-h)}^\alpha \Gamma(L^{h}\Theta)}& = 	\E[\big]{ \Gamma(  \Theta  )}, \quad h\in \Z^d }
for any $0$-homogeneous integrable functional $\Gamma: [0,\IF)^\Z \to  \R$ that vanished for $x \in [0,\IF)^\Z, x_0=0$. 
  If $\eZ(t),t\in \R^d$ is a max-stable stationary process with  marginals $\Phi_\alpha(x)= e^{-1/x^\alpha},x>0$, then clearly $(\eZ(t_1) \ldot \eZ(t_n)), t_i \in \Z^d, i\le n$ are for all $n\in \N$  regularly varying with index $\alpha$, and therefore $\eZ(t),t\in \Z^d$ has a STP which we denote by $\Theta$.  Below we specify  $\Theta$ in terms of $Z$ utilising i)  tilting (change of measure) for $V_\oldemptyset \fdd Z \lvert (Z(0)>0)$ and ii) the  {\it tilt-shift formula \eqref{pla}}. 

\BT \label{minor2} If $\eZ(t),t\in \R^d$ is a stationary max-stable process with unit marginals $\Phi_\alpha, \alpha>0$, then $\eZ$  has STP  $ \Theta(t)=\Xi_\oldemptyset V_\oldemptyset(t)\ge 0,t\in \Z^d$ and for any 0-homogeneous functional  $\Gamma:$  $\mathcal{F}_0\to \R$ which is $\mathcal{B}(\mathcal{F}_0) / \mathcal{B}(\R)$ measurable for any $h\in \R^d$ we have 
\bqn{ \label{pla} \quad 
	\E*{ \mathbb{I}\{Z(-h)>0\} Z^\alpha(-h) \Gamma(L^{h}Z)}
	= \E*{ Z^\alpha(\oldemptyset) \Gamma( Z)}=	\E[\big]{ \Gamma(  \Theta  )},
}
provided that the expectations exist.\\
\ET
By  the above clearly the SPT of $\eZ$ is non-negative, and 
if $\pk{Z(h)>0}=1,h\in \Z^d$, then  
$$\Theta\fdd Z.$$ 
Consequently, \eqref{pla} reduces to \eqref{tcf}. Note that $\Theta(t)= \Xi_\oldemptyset V_\oldemptyset(t)$ is defined for any $t\in \R^d$ and satisfies \eqref{pla} for all $t\in \R^d$, whereas the  {\it time-change formula} is stated in the literature only for discrete stationary time series. \\
 Conversely, if $\Theta(t) \in \TT$ is strictly positive with $\Theta(\oldemptyset)=1$ a.s.\ satisfying \eqref{tcf} for any $h\in \TT$, then by \netheo{ThMB} $\eZ(t),t\in \TT$ with $Z= \Theta$ is max-stable and stationary with STP $\Theta$.  \\

\COM{

Next, consider for simplicity $d=1$. 	  It is well-knwon that STP's are important for the calculation of extremal indices, which in case of a jointly  regularly varying $X(t),t\in \Z$ with index $\alpha>0$ is  equal to
 	 $$ \nu= \pk[\Big]{ \sup_{t \le -1, t\in \Z}  R \abs{\Theta(t)} \le 1}>0,$$ 
 	  with $R$ a unit Pareto rv independent of the STP $\Theta$ of $X$, provided that the finite mean cluster size (FMCS) condition \cite{BojanS}[Cond.\ 4.1] and the anti-clustering condition hold, see \cite{BojanS}. By 
 	  \cite{BojanS}[Prop.\ 4.2] the FMCS condition  implies 
 	  \bqn{ \label{m3T}
 	  	\pk[\Big]{ \limit{t\in \Z, \abs{t}  } {\Theta(t)}=0}=1.
 	  }
Let $\eZ(t), t\in \Z$ be a stationary max-stable time series  with marginals $\Phi_\alpha$ and STP $\abs{\Theta}$. Along the arguments of the proof of  \cite{KE17}[Thm. 2.1] we have that \eqref{m3T} implies that the FMCS condition holds for $\eZ$. Moreover, by \netheo{mentalist}, \netheo{minor2} and  \cite{BojanS}  the Pickands constant $\mathcal{H}_Z^1$ is equal to  $\nu$ and   
 	 	\BQN\label{Bi}
 	 	\nu = \mathcal{H}_Z^1= \E[\Big]{\frac{\sup_{t\in \Z} \abs{\Theta(t)}^\alpha} {\sum_{t\in  \Z} \abs{\Theta(t)}^\alpha}  }>0,
 	 	\EQN
provided that \eqref{m3T} is valid.\\
 To this end, let  $X(t),t\in \Z$ be a jointly regularly varying time series with SPT $\Theta$ and marginals $\Phi_\alpha$. In view of \cite{Janssen} and our findings, under \eqref{m3T} we can construct a max-stable stationary time series $\eZ(t),t \in \Z$ with SPT $\abs{\Theta}$ for which \eqref{Bi} holds.  \\
Therefore,  we conclude that formulas for  $\nu$ of a jointly regularly varying time series can be used to calculated Pickands constants of corresponding max-stable stationary processes with \FRE marginals. Moreover, the extremal indices of jointly regularly varying time series can be studied by restricting to the class of max-stable stationary time series. On the other hand, since Pickands constants are also defined for max-stable processes in continuous times, they provide a natural extension of the definition of extremal index of continuous time stationary processes, recall also Remark \ref{D_remDelta0}. 
	
}

\COM{
\subsection{Two-sided extensions of general max-stable processes} 
Let $\eZ(t), t\ge 0$ be a max-stable processes as in the Introduction. If 
$\eZ$ is stationary, then we can define  a max-stable process $\zeta_Y$ in the real line which has the same fidi's as $\eZ(t), t \ge 0 $, see \cite{MR3561100}[Prop. 1.1.2]. In other words, we can define a two-sided extension of $\eZ$ since we assume its stationarity. In Section 3, under the assumption that $Z(h)$ is a.s. finite $Y$ is determined explicitly such that its associated max-stable process $\zeta_Y$ is stationary and agrees with $\eZ$ for $t \in [0,\IF)$. Since $Z(h)$ is almost surely finite for some $h\ge 0$, and we can assume without loss of generality $\pk{Z(0)=0}=1$, alternatively the fidi's of $Y$ can be determined directly by $\zeta_Y$ using \netheo{Th1}. \\
If $\pk{Z(h)=- \IF}>0$, then only the identifiable part of $Y$ can be determined (using for instance \nelem{prop2}) i.e., we can determine $\ThetaH$. 
It is possible for this case to determine $Y$ as in \netheo{kolmog} by using  \eqref{kolmogorov} and the definition of $\Xi_h Z$ given in \eqref{spectWroc}. 
}
 	
 	\def\Z{\mathbb{Z}}
 	\def\inn{\in \mathbb{N}}
 	
 	{

 		\section{Proofs} \label{SecProof}

 		\prooflem{lem00}
 		For any $k>1, t_i \in \TT, x_i\inr, i \le n, x_{n+1}> 0$, the assumption $Y (h)=0$ almost surely implies (set $t_{n+1}=h, Y_i=Y(t_i), K=\{1 \ldot n+1\}$)
 		\BQNY
 		\lefteqn{\pk[\Big]{\forall j \in K:  \zeta_{Y}(t_j) \le x_j + \ln k \Bigl \lvert \zeta_{Y}(h)  > \ln k}}\\
&=&
 		\frac{  1  }{ 1- e^{- 1/k}}
 		\Bigl [
 		e^{-  \int_{\R}  \pk{ \exists j\in K : Y_j > x_j + \ln k -y}  e^{-y} \, dy }  -
 		e^{- \int_{\R}  \pk{ Y_{n+1} > \ln k- y, \text{ or }  \exists j \in K : Y_j > x_j + \ln k -y}  e^{-y} \, dy }\Bigr] \\
 		&=&
 		\frac{  1  }{ 1- e^{- 1/k}}
 		\Bigl[
 		e^{-  \frac{1}{k} \int_{\R}  \pk{ \exists j\in K: Y_j > x_j -y}  e^{-y} \, dy }  -
 		e^{- \frac{1}{k}\int_{\R}  \pk{ Y_{n+1} > - y, \text{ or }  \exists j\in K: Y_j > x_j  -y}  e^{-y} \, dy }\Bigr] \\
 		& \to &  \int_{\R} \Bigl[
 		\pk{y > 0, \text{ or }  \exists  j \in K : Y_j > x_j  -y} -
 		\pk{ \exists j\in K  : Y_j > x_j  -y}\Bigr]  e^{-y} \, dy,  \quad k\to \IF\\
 		& = &  \int_{\R} \Bigl[  \pk{ \forall j\in K : Y_j \le   x_j  -y}- \pk{y \le  0, \forall j\in K : Y_j \le  x_j  -y}
 		\Bigr]  e^{-y} \, dy\\
 		&=&
 		\int_{0}^\IF \pk{\forall  j\in K : Y_j \le  x_j  -y}  e^{-y} \, dy \\
 		&= &\pk{ \forall j \in K: Y_j+ \mathcal{E}  \le x_j},
 		\EQNY
 		where $\mathcal{E}$ has a unit exponential df being independent of $Y$. Consequently, since a.s. $Y_h=Y(h)=0$ we have the convergence in distribution as $k\to \IF$
 		\bqn{ \Bigl(\zeta_Y(t_1)- \zeta_Y(h) \ldot \zeta_Y(t_n)- \zeta_Y(h), \zeta_Y(h) \Bigr) \Bigl \lvert (\zeta_Y(h)> \ln k)
 			\todis (Y_1 \ldot Y_n, \mathcal{E}),  \label{lC}
 		}
 		hence the proof is complete.  \QED
 		
 		\prooftheo{Th1} $i)$ Let $H$ be the df of $(\eZ(t_1) \ldot \eZ(t_n)) $. For ${W}$ with df  $G$ given in \eqref{dfG} set
$(W_{1}^{(h)}\ldot W_{n}^{(h)})= {W} \lvert (W_h>0).$ 
If  $h\in \{t_1 \ldot t_n\}$, then  by \cite{Segers}[Eq. (8.67)]  
$$ \Bigl( \eZ(t_1) - \ln k \ldot \eZ(t_n)  - \ln k \Bigr) \Bigl \lvert \bigl ( \eZ(h) > \ln k \bigr) \todis (W_1^{(h)}  \ldot W_n^{(h)}) , \quad k\to \IF$$
holds,  see also   \cite{raey}[Eq. (13)]. (Note that $W_k^{(h)}, k\not=h$ may assume value $-\IF$ with non-zero probability).  Hence the claim follows by \eqref{lC}. \\
$ii)$ First note that $\E{e^{Z(t)}}=1, t\in \TT$
 		implies for $\mu$ a probability measure on $\TT$ (recall $\SZ=  \int_{\TT} e^{Z(s)} \mu(ds)$)
 		$$ \E{ \SZ} = \int_{\TT} \E{ e^{Z(s)} }\mu(ds)=1$$
and a.s. $\SZ < \IF$, which is assumed to hold if    $\mu$ is a positive  $\sigma$-finite measure. Since a.s.\ $Z(h)> - \IF ,h\in \TT$, then  a.s.\ $\SZ >0 $. Hence,  for any $h, t_i\in \TT, x_i \inr, i\le n$  by 
Fubini-Tonelli theorem
 		\bqn{ \label{les}
  			- \ln H(x)
 			&=& \E[\Big]{\frac{ \SZ}{\SZ} \max_{1 \le k \le n} e^{- x_k+ Z(t_k) }} \notag\\
 			&=& \int_{\TT} \E* {  e^{Z(h)} \frac{ \max_{1 \le k \le n} e^{- x_k+ Z(t_k) - Z(h)  }}
 				{ \int_{\TT} e^{ Z(s)- Z(h)} \mu(ds) }} \, \mu(dh)\\
 			&=& \int_{\TT} \E*{  \max_{1 \le k \le n} e^{- x_k+ \Dhz(t_k)  - \ln \bigl(\int_{\TT} e^{ \Dhz(s)} \mu(ds)\bigr) }} \, \mu(dh). \notag
  		}
Now, if $\sum_{i\ge 1} \varepsilon_{(P_i, T_i, Z_i)}$ is a PPP on $\R \times \TT \times \R^\TT$ with intensity $e^{-p} dp \cdot  \mu(dt) M(dz)$ where
$M$ is the law of $Z$, then for
$$\eta(t)=\max_{i\ge 1} \Biggl( P_i   -\ln \int_{\mathcal{T}} e^{  \Xi_{T_i} Z_i(s)  - \Xi_{T_i} Z_i(t)} \, \mu(ds) \Biggr)$$
we have
\bqny{- \ln \pk{\eta(t_i) \le x_i, 1 \le i \le n} &=&  \E*{\int_{\TT}   \max_{1 \le k \le n} e^{- x_k+ \Dhz(t_k)  - \ln \bigl(\int_{\TT} e^{ \Dhz(s)} \mu(ds)\bigr) } \, \mu(dh)},
	}
hence the claim follows.
 		\QED

 	\def\inn{\in \mathbb{N}}
 	 
 	\def\hh{ \thhh }
 	
 	\def\TT{\mathcal{T}}

 	\prooftheo{ThM}  \underline{$a)$ implies $d)$}:\\
 	Let $t_i\in \TT,i\le n $ be given. As mentioned in the Introduction for the validity of \eqref{repA}, since  a.s. both $Z(a)> -\IF$ and $Z(a+h)> - \IF$ hold, then we have 
 	$$\eZ
 	\fdd   \zeta_{ \Xi_{a} Z} \fdd   \zeta_{ L^h\Xi_{a} Z}   \fdd  \zeta_{ \Xi_{a+h} Z}.$$  
	By our definition in \eqref{abu}   
 	$$
  L^h\Xi_{a}	 Z (t)=  L^h (  Z^{[a]}(t)- Z^{[a]}(a))=  Z^{[a]}(t-h)- Z^{[a]}(a), \quad t\in \TT
 	$$
we have a.s.
 	$
 	 L^h \Xi_a   Z(a+h)= \Xi_{a+h} Z(a+h)=0$. Hence  by \netheo{Th1} $Y\EQD Y^*$.\\
 	\underline{$d)$ implies $c)$}:\\
 	First note that by the shift-invariance of $\Gamma$
 	\bqny{
 		\E*{\Gamma( \Xi_{a+h} Z)}&=& \E*{e^{Z(a+h)} \Gamma\bigl(Z- Z(a+h)\bigr)}=\E*{ e^{Z(a+h)} \Gamma(Z)}
 	}
 	and
 	\bqny{
 		\E*{\Gamma( L^h\Xi_{a} Z )}&=& \E*{e^{Z(a)} \Gamma\bigl( \thhh  Z-  Z(a)\bigr)}
 		= \E*{e^{Z(a)} \Gamma( \thhh  Z)}. 
 	}
 	By  statement $ d)$ we have that the fidi's  of
 	$ \Xi_{a+h} Z$ and $ L^h\Xi_a Z$ are the same, which together with the measurability of  $\Gamma$ implies
 	$ \Gamma\bigl( \Xi_{a+h} Z\bigr) \equaldis \Gamma\bigl(  L^h\Xi_{a}  Z\bigr).$
 	 Consequently, we obtain
 	\bqny{
 		\E*{\Gamma\bigl( \Xi_{a+h} Z\bigr)}&=&
 		\E*{\Gamma\bigl( L^h\Xi_{a}  Z)\bigr)}=\E*{ e^{Z(a+h)} \Gamma(Z)}=
 		\E*{e^{Z(a)} \Gamma( \thhh  Z)}
 	}
 	and thus the claim follows. \\
 	\underline{$c)$ implies $b)$}:
 	
If $\mu$ is a probability measure on $\TT$,
 	for any $h,t_i \in \TT, x_i\inr, i\le n$  by \eqref{thmd},  Fubini-Tonelli theorem, statement $c)$ and \eqref{les} yield
 	\bqny{
 		- \ln \pk{ \eZ(t_i) \le x_i, i= 1\ldot n}
 		&=& \int_{\TT} \E*{  e^{Z(h)} \frac{ \max_{1 \le k \le n} e^{- x_k+ Z(t_k)  }}
 			{ \int_{\TT} e^{ Z(s)} \mu(ds) }} \, \mu(dh)\\
 		&=:& \int_{\TT} \E*{  e^{Z(h)} \Gamma(Z)}\, \mu(dh)\\
 		&=& \int_{\TT} \E*{  \Gamma(  L^h\Xi_{\oldemptyset} Z  ) }\, \mu(dh)\\
 		&=& \E*{  \Gamma(  L^T \Xi_{\oldemptyset}  Z  ) },
 	}
 	with $T$ a copy of $T_1$, which is independent of $\Xi_{\oldemptyset} Z$ (recall  $\oldemptyset= (0 \ldot 0)$). Consequently,
 	\bqny{
 		- \ln \pk{ \eZ(t_i) \le x_i, i= 1\ldot n}
 		&=& \E*{ \max_{1 \le k \le n} e^{-x_k + L^T \Xi_{\oldemptyset}  Z(t_k)    -
 				\ln \int_{\TT} e^{  L^T\Xi_{\oldemptyset} Z (s)} \mu(ds)} }.
 	}
The case that $\mu$ is  a positive $\sigma$-finite measure on $\TT$ such that \eqref{trift} holds follows with similar arguments. \\
 	\underline{b) implies a)}:\\ Let $\mu$ be the Dirac measure at $h$. We have
 	that $\eZ$ has the same law as $\max_{i\ge 1}(P_i+  \Xi_{\oldemptyset} Z_i(t-h)),t\in \TT $ and by \eqref{repA}, this implies that
 	$\eZ \fdd \zeta_{ \thhh  Z}$, hence $\eZ$ is stationary. \\
We note that applying statement $c)$ to the functional 
 	$\Gamma_A(Z)=\mathbb{I}\{ (W(t_1)\ldot W(t_n))\in A\}$ with 
$W(t_i)=Z(t_i)- Z(a+h)$ yields 
\bqny{\pk{( \Xi_{a+h} Z(t_1) \ldot  \Xi_{a+h} Z(t_n)) \in A  } &= &
 		\E*{e^{Z(a+h)} \Gamma_A(Z) } =  		\E*{e^{Z(a)} \Gamma_A( \thhh  Z) }  \\
 		& = &
\pk{(  L^h\Xi_{a} Z(t_1) \ldot  L^h\Xi_{a}Z(t_n)\in A  },
 	}
which implies statement $d)$.
 	\QED
 	
\prooftheo{korrK} Let $\varphi(t)= \ln \E{e^{(t,X)}},t\inr^d$ denote  the cumulant generating function of $X$. As shown in   \cite{MR3454027}[Th.\ 1] $\eZ$ is stationary
in $\R^d$ i.e., in our setup $Z(t)= (t, X)  - \kappa(t),t \inr^d$  is Brown-Resnick
stationary, if and only if $X$ is Gaussian with mean $\mu $, covariance matrix $\Sigma$ and further $\kappa(t)= \varphi(t),t\inr^d$.
Our assumption is slightly weaker since we assume the stationarity of $\eZ$ restricted on $\Zdd$ for any $\delta>0$ and not its stationarity on $\R^d$.\\
 As in the proof of Theorem 1 therein, our assumption is equivalent with
\bqn{ \label{lolby}
	\varphi\Bigl(\sum_{i=1}^n u_i t_i\Bigr)- \sum_{i=1}^n u_i \varphi(t_i) = \varphi\Bigl(h + \sum_{i=1}^n u_i t_i \Bigr)- \sum_{i=1}^n u_i \varphi(h+t_i)}
for any $h,t_i \in \Zdd,u_i \inr, i\le n$ where $\sum_{i=1}^n u_i=1$. Write next $[a]=([a_1] \ldot [a_d])$ for any $a\inr^d$ with $[a_i]$ the largest integer smaller than $a_i$.  By \eqref{lolby} for any $\lambda\in (0,1),\delta>0$  and any $t_0, t_1,t_2 \in \R^d$  (set $h= \delta t_0,
z_i= \delta [ t_i/\delta],   v= \lambda  z_1  + (1- \lambda )
z_2  $)
\bqny{ \varphi( v )- \lambda
	\varphi( z_1 )-
(1- \lambda ) \varphi( z_2  )
& =&
	\varphi( h +  v )- \lambda
	\varphi(  h+ z_1)-
	(1- \lambda ) \varphi(h+ 	z_2  ) .
		}
 Letting $\delta \to 0 $ we obtain
$$ \nabla \varphi \bigl( \lambda t_1 + (1- \lambda ) t_2 \bigr)= \lambda \nabla \varphi(t_1) +
(1- \lambda ) \nabla \varphi(t_2),
 $$
where $\nabla \varphi(t)$ denotes the gradient of $\varphi$ at $t\inr^d$. Consequently, with the same arguments as in \cite{MR3454027}  it follows that $X$ is Gaussian with mean $\mu $ and covariance matrix $\Sigma$ and further $\kappa(t)= \varphi(t),t\inr^d$. Hence the proof follows by Remark \ref{remiii} $ii)$.
\QED

\prooftheo{korrP} If $B$ has stationary increments, then  by \eqref{dusig}
\bqn{
	\label{sorgmach}
	\Dhz \fdd  \thhh  Z, \quad \forall h \in \TTT,
}
hence statement $d)$ in \netheo{ThM} implies that $\eZ$ is stationary.
Note that \eqref{sorgmach} is previously shown in 
\cite[Prop. 2]{kab2009a}. Conversely, if $\eZ$ is stationary, then by statement $d)$ in \netheo{ThM} we have that \eqref{sorgmach} holds, which combined with \eqref{exDhZ} yields for any $h\in \TT$
\bqny{
	B(t-h) - \sigma^2(t-h)/2\fdd B(t)- B(h) - Var(B(t)- B(h))/2, \quad t\in \TT,
}
hence $B$ has stationary increments and thus the claim follows.
\QED

 	\prooftheo{kolmog}   The proof is based on the result of \netheo{ThM}, which is also valid for  $\TT= [0,\IF)$.
 First, we show that $Y$ can be defined using the  Kolmogorov's consistency theorem, see e.g., \cite{YuliaK}[Th.\ 1.1].  	It suffices to consider in the following only $t_1< \cdots < t_n\inr$ such that $t_1 < 0$.
For any permutation $\pi$ of
 	$t_1 \ldot t_n$ we have that
 	$$ Y_{\pi(t_1) \ldot \pi(t_n)}= \Bigl(Y(\pi(t_1))  \ldot Y(\pi(t_n))\Bigr) \EQD
 	\Bigl( \Xi_{h}(\pi(t_1))   \ldot  \Xi_{h}(\pi(t_n))  \Bigr), \quad h=-  t_1$$
 	implying  $Y_{\pi(t_1) \ldot \pi(t_n)} \EQD Y_{t_1 \ldot t_n}$ since $h$ is independent of the chosen permutation. \\
 	The consistency of the family of fidi's  follows if we can further show that
 	for any non-empty $I \subset \{1 \ldot n\}$ (write $J:= \{1 \ldot n\}\setminus I$)
 	\BQNY
 	\pk{ Y(t_i) \le x_i, i\in I, Y(t_j) \inr, j\in J}&=& \pk{ Y(t_i) \le x_i, i\in I}
 	\EQNY
 	for any $x_i\inr, i\in I$.  
 	If $1 \in I$, the above follows immediately by the definition of $Y$.  Suppose next that $1\in J$ and assume for notation simplicity that $J=\{1 \}$.
 	We need to show that
 	\BQNY
 	\pk{ Y(t_1)\inr, Y(t_i) \le x_i, 2 \le i\le n}&=&
 	\pk{L^{t_1} \Xi_{-t_1}  Z(t_i)  \le x_i, 2 \le i\le n}\\
 	&=&    \pk{ L^{t_2} \Xi_{-t_2} Z(t_i) \le x_i,  2 \le i\le n}\\
 	&=&  	\pk{Y(t_i) \le x_i, 2 \le i\le n}
 	\EQNY
 	for any $x_i\inr, 2 \le i\le n$, which follows directly by 
 	\COM{  Indeed, this is the case since if $t_2 < 0$, then
 	\BQNY
 	\Bigl(Y(t_2) \ldot Y(t_n)\Bigr) &\EQD &
 	\Bigl( \Xi_{-t_2} Z(0) \ldot  \Xi_{-t_2} Z(t_n- t_2) \Bigr)\\
 	&\EQD  &   \Bigl(   L^{-(t_2-t_1)}\Xi_{-t_1}  Z(0) \ldot  L^{-(t_2-t_1)}\Xi_{-t_1} 
 	Z(t_n- t_2)\Bigr)\\
	&\EQD &   \Bigl(  \Xi_{-t_1}  Z(t_2-t_1) \ldot  \Xi_{-t_1}  Z(t_n-t_1) \Bigr),
 	\EQNY
 	where the second equality in distribution follows from}  \eqref{thmd}. 
 \COM{Similarly, for $t_2>0$ since $h=0$, then
 	\BQNY
 	\Bigl(Y(t_2) \ldot Y(t_n)\Bigr)
 	&\EQD &    \Bigl( \Xi_{0} Z(t_2) \ldot  \Xi_{0} Z(t_n)\Bigr)\\
 	&\EQD &   \Bigl(  L^{t_1}\Xi_{-t_1}   Z(t_2) \ldot  L^{t_1}\Xi_{-t_1}  
 	Z(t_n)\Bigr)\\
&\EQD &   \Bigl(  \Xi_{-t_1}  Z(t_2-t_1) \ldot  \Xi_{-t_1}  Z(t_n-t_1) \Bigr),
 	\EQNY
 }
\COM{ 
 \BQNY 
\arccos(-1+ \ve)- \arccos(-1)&=&  g(\ve)- g(0) \\
&=& \int_0^\ve g'(z) dz\\
&=& \int_0^\ve - (1-  (1- x)^2 )^{1/2- 1}\ dx\\
&=& -\int_0^{\ve} x^{1/2- 1} ( 2  - x)^{1/2- 1}  dx\\
& \sim &  - \frac{1}{\sqrt{2}} \int_0^{\ve} x^{1/2- 1}    dx \\
&  \sim&  -   \sqrt{2 \ve}, \quad \ve \to 0
 \EQNY	
} 	
Hence since the conditions of
 	\cite{YuliaK}[Th.\ 1.1] are satisfied, then $Y(t),t\inr$ exists. 
 By \eqref{kolmogorov} for any $t\inr$
 	$$ Y(t)\EQD  \Xi_{-t} Z(0)= Z^{[-t]}(0)- Z^{[-t]}(-t)$$
 	implying that $ 	\E[\big]{e^{Y(t)}}= 1$, 
\COM{ 	\BQNY
 	\E[\big]{e^{Y(t)}}= \E[\big]{e^{Z^{[-t]}(0)- Z^{[-t]}(-t)}}= \E[\big]{e^{Z(-t)}e^{Z(0)- Z(-t)}}=
 	\E[\big]{e^{Z(0)}}=1,
 	\EQNY
}
hence  $\zeta_{Y}(t),t\inr $ associated to $Y$ (as in \eqref{W} with $Z$ substituted by $Y$ and $\mathcal{T}=\R$) has
 unit Gumbel marginals and is max-stable.  The stationarity of 
 $\zeta_{Y}$ follows easily, we omit the proof. 
 \COM{arg1is stationary. 	For any $t_1 \ldot t_n \inr, -h_1:=t_1 \le t_j,1 < j\le n$  by the definition 
 	\BQNY
 	\Bigl (Y(t_1) \ldot Y(t_n)\Bigr)
 	&\EQD&
 	\Bigl ( \Xi_{h} Z(0) \ldot  \Xi_{h} Z(t_n- t_1)\Bigr) .
 	\EQNY
 	Hence since 
 	$\eZ \fdd  \zeta_{\Dhz }$, then
 	$$ \Bigl( \zeta_{Y}(t_1) \ldot \zeta_{Y}(t_n)\Bigr) \EQD
 	\Bigl( \zeta_{Z}(0) \ldot \zeta_{Z}(t_n-t_1)\Bigr).$$
 	Further
 	 	\BQN \label{ang1}
 	\Bigl (Y(0) \ldot Y(t_n-t_1)\Bigr)
 	&\EQD &
 	\Bigl ( \Xi_{0} Z(0) \ldot  \Xi_{0} Z(t_n-t_1)\Bigr)
 	\EQN
 	yields 
\BQN\label{ang2}
\Bigl( \zeta_{Y}(0) \ldot \zeta_{Y}(t_n-t_1)\Bigr) \EQD
 	\Bigl( \zeta_{Z}(0) \ldot \zeta_{Z}(t_n-t_1)\Bigr).
 	\EQN
Consequently, 
by \eqref{ang1} and \eqref{ang2},  for any $a\inr$
 	
 	$$ \Bigl( \zeta_{Y}(t_1+a) \ldot \zeta_{Y}(t_n+a)\Bigr) \EQD
 	\Bigl( \zeta_{Y}(0) \ldot \zeta_{Y}(t_n-t_1) \Bigr)
 	\EQD  \Bigl( \zeta_{Y}(t_1) \ldot \zeta_{Y}(t_n)\Bigr) $$
 	and thus $\zeta_{Y}$ is max-stable and stationary establishing the proof.
 	} \QED 

 	\prooflem{lemDhV}
 	For notational simplicity write $J,V,W$ instead of
 	$J_h, V_h,W_h$.  Since $J $ is independent of $V$ and a.s. $W(h)=-\IF$, the assumption that $\E{e^{Z(h)}}=1$  	implies (recall $0 \cdot \IF=0$ and set $p= \pk{J=1}$)
\bqn{ \label{nach} 1= \E[\big]{e^{Z(h)}}=\E[\big]{ e^{J V(h)+(1-J) W(h)  }} =
 	\E[\big]{ \mathbb{I}\{ J=1\} e^{ J V(h) }}=\E[\big]{ e^{ V(h) + \ln p}}.}
 	For any $t_1=h \ldot t_n\in \mathcal{T},x\inr^n $ we have (set $\bar J= 1- J$)
 	\BQNY
 	- \ln \Hxn 
 	&=& \E[\Big]{   \max_{1 \le j \le n} e^{J V(t_j)+\bar J W(t_j) - x_j }  }\\
 	&=& \E[\Big]{  e^{ V(h) + \ln p}  \max_{1 \le j \le n} e^{J [V(t_j)  - V(h)]+\bar J [W(t_j) - V(h)]- \ln p - x_j}  }\\
	&=& 
 	\E*{  \max_{1 \le j \le n} e^{J [V ^\HH (t_j)  - V ^\HH (h)]+\bar J [W(t_j) - V ^\HH (h)]- \ln p - x_j}  }, 
 	\EQNY
where we used the fact that  both $J,W$ are independent of $V$. Hence we have  $\eZ \fdd  \zeta_{ \Dhz}$.\\
 	 Next, consider the PPP $ \sum_{i\in N}\varepsilon_{(P_i, T_i, Z_i)}$ on $\R \times \TT \times  \R^{\TT}$.
 	Using the void probability formula and $\eZ \fdd  \zeta_{ \Dhz}$,  for any $t_i\in \TT, x_i \inr, i\le n $  we obtain
 	\bqny{- \ln \pk{ \eta(t_i) \le x_i, 1 \le i \le n}&=&
 		\int_{\TT} - \ln \pk{ \zeta_{\Dhz}(t_i) \le x_i, 1 \le i \le n}
\mu(dh)\\
 		&=&   - \int_{ \TT} \ln \pk{ \eZ(t_i) \le x_i, 1 \le i \le n}  \mu(dh)\\
 		&=&   - \ln \pk{ \eZ(t_i) \le x_i, 1 \le i \le n} ,
 	}
establishing the proof. \QED

 	\prooftheo{Euler} Define next for a process $X(s),s\in \TT$ 
 	$$ Q_{t, x}(X)= \inf {\rm argmax}_{ 1  \le  j \le n} e^{ X(t_j)  - x_j}, \quad t=(t_1 \ldot t_n) \in \TT^n,\quad  x  \inr^n.$$
 	Hereafter write  $\E{K; B}:=\E{K\mathbb{I}\{B\}}$ for $K$ some random element and $B$ an event. Recall that $W_h(h)=- \IF$ a.s. and  set
\bqn{ \label{yh}  Y_h(t):=J_hV_h(t) + (1-J_h)W_h(t) \fdd  Z(t), \quad h\in \TT.}
For any $x\inr^n$ we have (recall that a.s. $W_h(h)=-\IF$ and the indicator rv's $J_{t_k}$ are independent of $V_{t_k},W_{t_k}$)
 	\BQN
\lefteqn{ 	- \ln H(x)} \notag\\
 	&=&\E*{ \sum_{k=1}^n \mathbb{I}\{ Q_{t, x}(Z)  =k\}
 		\max_{ 1\le j \le  n} e^{ Z(t_j)- x_j} }
 	\notag \\
 	&=&\sum_{k=1}^n e^{- x_k} \Bigl[ \E[\Big]{ \mathbb{I}\{ Q_{t,x}(Y_{t_k}) =k \} e^{ Y_{t_k}(t_k) }; J_{t_k}=1 }
 	+ \E[\Big]{ \mathbb{I}\{  Q_{ t,x}(Y_{t_k}) =k \} e^{Y_{t_k}(t_k)}; J_{t_k}=0 }\Bigr]
 	\notag \\
 	&=&\sum_{k=1}^n \pk{J_{t_k}=1} \E[\Big]{ \mathbb{I} \{ Q_{t,x}(V_{t_k}) =k \}
 		e^{ V_{t_k}(t_k) }} \notag \\
 	&=&\sum_{k=1}^n e^{- x_k}\E[\Big]{ \mathbb{I}\{ Q_{t,x}(V_{t_k})=k \}
 		e^{ V_{t_k}(t_k) + \ln \pk{J_{t_k}=1} }} \notag \\
 	&=&\sum_{k=1}^n e^{- x_k}\E[Big]{ \mathbb{I}\{ Q_{t,x}(V_{t_k}^{[t_k]})=k \} }  \label{tauVh}\\
 	&=&\sum_{k=1}^n e^{- x_k} \pk[\Big]{ \max_{1 \le i < k} \Bigl( V^{[t_k]}_{t_k}(t_i)- x_i\Bigr)< V^{[t_k]}_{t_k}(t_k)- x_k,
 		\max_{ n\ge  i >     k} \Bigl(V^{[t_k]}_{t_k}(t_i)- x_i\Bigr) \le V_{t_k}^{[t_k]}(t_k)- x_k} \notag \\
 	 	&=&\sum_{k=1}^n e^{- x_k} \pk[\Big]{ \inf {\rm argmax}_{1 \le i \le n} \Bigl(  \Theta_{t_k} (t_i)  -  x_i\Bigr) = k} , \notag
 	\EQN
where  $ \ThetaH(t):= \Xi_h V_h= V_h ^\HH (t)- V ^\HH _h(h)$ and $\max_{1 \le j< 1}(\cdot)=\max_{n \ge i> n}(\cdot) =: -\IF$, 	hence the claim follows.  \QED

 	\prooftheo{ThMB}   \underline{$a)$ implies $c)$}:\\
Let $\OO Z$ be given from \eqref{spectWroc}. 
As in the proof of \netheo{ThM} for any $a,h\in \TT$ we have $ \zeta_{  \Xi_{a+h} Z} \fdd  \zeta_{ L^h \Xi_a    Z},$ hence   by \nelem{prop2} in Appendix  
 	$$  \Theta_{a+h}:=\Xi_{a+h} V_{a+h} \fdd  L^h\Xi_{a}   V_a=: L^h \Theta_a  $$  establishing the claim.\\
 	\underline{$c)$ implies $b)$}:\\
 	 First note that by the shift-invariance of $\Gamma$, for any $h\in \TT$ we have
 	$$\Gamma(V_h ^\HH )=\Gamma( V_h ^\HH - V ^\HH (h))= \Gamma( \ThetaH).$$
Further, since $J_h$ is independent of $V_h$ and $W_h$, (recall $W_h(h)=-\IF$ a.s.), then using the shift-invariance of $\Gamma$ yields
 	\bqny{
 		\E[\big]{e^{Z(h)}\Gamma( Z)}&=&  \E[\big]{ e^{J_hV_h(h)+ (1- J_h) W_h(h)} \Gamma(J_hV_h+ (1- J_h) W_h)}\\
 		&=& \E[\big]{ e^{V_h(h) + \ln \pk{J_h=1}} \Gamma(V_h)}=
 		\E{ \Gamma(V_h ^\HH )}= \E{ \Gamma( \ThetaH)}.
 	}
 	Consequently, since $ \ThetaAH \fdd   L^h \Theta_ a $ is valid
 	for any $a,a+h\in \TT$, then
 	\bqny{
 		\E[\big]{e^{Z(a+h)}\Gamma( Z)}&= \E[\big]{ \Gamma( \Theta_{a+h} )}=
 		\E[\big]{ \Gamma( L^h\Theta_{a} )} =  \E[\big]{e^{Z(a)}\Gamma(  \thhh  Z) }
 	}
 	establishing the claim. \\
 	\underline{$b)$ implies $a)$}:\\
 	  Given $t_1 \ldot t_n \in \TT$ distinct and $h\in \TT$ by  \eqref{tauVh}
 	  for any $x \inr^n$ 
 	\bqny{
 		- \ln \Hxn &=& -\ln \pk{\eZ(t_i) \le x_i, 1 \le i\le n} =  \sum_{k=1}^n e^{- x_k}\E[\big]{\Gamma_k( \Theta_{t_k} )}, 
 	}
 	where 
 	$$ \Gamma_k( \Theta_{t_k}):=
 	\mathbb{I}
 	\bigl\{
 	\inf {\rm argmax}_{1 \le i \le n}
 	\bigl(  \Theta_{t_k}(t_i)  - x_i\bigr) = k \bigr\}. $$
 	The functional $\Gamma_k$ is shift-invariant, hence from statement $b)$
 	$$ \E{\Gamma_k(  \Theta_{t_k})}= \E[\big]{\Gamma_k(  L^h\Theta_{t_k-h}  )} $$
 	implying (set $t_i^*= t_i- h$)
 	
 	\bqny{
 		- \ln \Hxn 
 		&=& \sum_{k=1}^n e^{- x_k}\E[\big]{\Gamma_k( L^h\Theta_{t_k-h} )}\\
 		&=& \sum_{k=1}^n e^{- x_k}\E[\Big]{ \mathbb{I}\bigl\{
 			\inf {\rm argmax}_{1 \le i \le n}
 			\bigl(  \Theta_{t_k-h}(t_i- h)  - x_i\bigr) = k \bigr\}}\\
 		&=& \sum_{k=1}^n e^{- x_k}\E{ \Gamma_k(  \Theta_{t_k^*} )} = - \ln \pk{\eZ(t_i-h) \le x_i, 1 \le i\le n},
 	}
 	which proves the stationarity of $\eZ$, hence the claim follows.  \QED

\COM{ 	\proofprop{th2} Since by the assumptions
 	$$ \limit{k} k \pk{X_i> a_k x+ b_k}= e^{-x}, \quad x\inr$$
 	we have for any $x_1 \ldot x_n \inr$ and $k\inn$
 	\BQNY
 	\lefteqn{k [1- F(a_kx_1 + b_k \ldot a_k x_n+b_k)]}\\
 	&=&
 	k\pk{ \cup_{1 \le j \le n} \{  X_j -a_k x_j> b_k \}}\\
 	&=&
 	k \pk{ \cup_{h=1}^n \{ \inf {\rm argmax}_ {1 \le j \le n}  ( X_j-a_k x_j)  = h\}, \cup_{1 \le j \le n} \{  X_j -a_k x_j> b_k \}}\\
 	&=&
 	k \sum_{h=1}^n \pk{\inf {\rm argmax}_ {1 \le j \le n}  ( X_j-a_k x_j)  = h, \cup_{1 \le j \le n} \{  X_j -a_k x_j> b_k \}}\\
 	&=&  \sum_{h=1}^n k \pk{ \inf {\rm argmax}_ {1 \le j \le n}  ( X_j-a_k x_j)  = h, X_h - a_k x_h > b_k }\\
 	&=&
 \sum_{h=1}^n  k \pk{X_h >  a_k x_h + b_k  }
 	\pk{ \inf {\rm argmax}_ {1 \le j \le n}  ( X_j-a_k x_j)  = h\Bigl \lvert  \frac{X_h - b_k}{a_k}> x_h }\\
 	&\to &  \sum_{h=1}^n  e^{- x_h} \Psi_h(x_1-x_h \ldot x_n- x_h), \quad k\to \IF,
 	\EQNY
 	hence the claim follows. \QED
}

\COMC{The next result extends Lemma 2 in \cite{SBK} formulated for the case $d=1$ under the assumption that a.s.   $Z(\oldemptyset)=0$.
	\COMB{
		\BEL \label{lemD} {Let $Z$ be as in \netheo{ThMB}. If $\mu$ is the  counting measure on  $(k\delta) \Zd,d\ge 1 $ with $k\inn, \delta > 0$,  then for any $T>0$}
		\BQN\label{cor2}
		\E[\Big]{ \sup_{t\in \delta \Zd \cap [0,T]^d} e^{ Z (t)}} 		&= &     T^d\int_{ [0,1]^d} \mathbb E\left\{\frac{\sup_{t\in \delta\Zd \cap \times_{i=1}^d[-h_i T,(1-h_i) T] } e^{  \Theta (t)}}
		{ \int_{  \delta\Zd \cap \times_{i=1}^d[-h_i T,(1-h_i) T] } e^{ \Theta(s)}\mu(ds) } \right\}\, \mu^T(dh) ,
		\EQN
		with $\Theta= \Xi_\oldemptyset V_\oldemptyset, \mu^T (dh)= \mu ( T dh)/T^d$ and $h=(h_1 \ldot h_d)$.
		\EEL
		\def\Ed{\mathcal{E}_\delta}
		\prooflem{lemD}
		Since $\mu$ is a counting measure on  $(k\delta) \Zd \cap [0,T]$, then (set $\Ed:=  \delta \Zd \cap [0,T]^d$)
		$$\pk*{\int_{ \Ed} e^{Z(s)} \mu (ds)=0,  \sup_{t\in \Ed} e^{Z(t)}>0  }=0.$$
		Consequently, with $Y_h$ defined in \eqref{yh}
		we have (recall $W_h(h)=-\IF$ a.s.)
		\bqny{
			\E[\Big]{ \sup_{t\in \Ed} e^{ Z (t)}}
			&= &  \int_{ [0,T]^d}
			\E*{ e^{Z(h)} \frac{ \sup_{t\in \Ed} e^{ Z (t)}}{\int_{ \Ed} e^{Z(s)} \mu (ds) } ;
				\sup_{t\in \Ed} e^{Z(t)}> 0   }\mu (dh)\\
			&= &  \int_{[0,T]^d}
			\E*{ e^{Y_h(h)} \frac{ \sup_{t\in \Ed} e^{ Y_h (t)}}{\int_{  \Ed} e^{Y_h(s)} \mu (ds) } ;
				\sup_{t\in \Ed} e^{Y_h(t)}> 0   }\mu (dh)\\
			&= &  \int_{[0,T]^d} \Biggl[
			\E*{ e^{Y_h(h)} \frac{ \sup_{t\in \Ed} e^{ Y_h (t)}}{\int_{  \Ed} e^{Y_h(s)} \mu (ds) } ;
				\sup_{t\in \Ed} e^{Y_h(t)}> 0; J_h=1   }\\
			&= &  \int_{ [0,T]^d}
			\E*{ e^{V_h(h)} \frac{ \sup_{t\in \Ed} e^{ V_h (t)- V_h(h)}}{\int_{  \Ed} e^{V_h(s)- V_h(h)} \mu (ds) } ;
				J_h=1   }\mu (dh)\\
			&= &  \int_{  [0,T]^d}
			\E*{ e^{V_h(h) - \ln \E{V_h(h)}} \frac{ \sup_{t\in \Ed} e^{ V_h (t)- V_h(h)}}{\int_{  \Ed} e^{V_h(s)- V_h(h)} \mu (ds) }
			}\mu (dh)\\
			&= &   \int_{   [0,T]^d}
			\E*{  \frac{ \sup_{t\in \Ed } e^{  \ThetaH (t) }}
				{\int_{ \Ed} e^{ \ThetaH(s)} \mu (ds) }    }\mu (dh)\\
			&= & \int_{  [0,T]^d}
			\E*{  \frac{ \sup_{t\in \Ed } e^{  \Theta (t-h) }}
				{\int_{  \Ed} e^{ \Theta(s-h)} \mu (ds) }    }\mu (dh),
		}
		where the last equality follows from $ \ThetaH \fdd  L^h \Theta
		 =L^h \Theta_\oldemptyset, $ which is a consequence of \netheo{ThMB}.
		Alternatively, using directly \eqref{xYx} and omitting few details, we obtain 
				\bqny{
			\E[\Big]{ \sup_{t\in \Ed} e^{ Z (t)}}
			&= &  \int_{ [0,T]^d}
			\E*{ e^{Z(h)} \frac{ \sup_{t\in \Ed} e^{ Z (t)}}{\int_{ \Ed} e^{Z(s)} \mu (ds) } ;
				\sup_{t\in \Ed} e^{Z(t)}> 0   } \mu (dh)\\
			&=:& 			\E*{ e^{Z(h)} \Gamma(Z)} =  \int_{  [0,T]^d} \E*{  \frac{ \sup_{t\in \Ed } e^{ L^h \Theta (t) }}	{ 
	\int_{  \Ed} e^{ L^h\Theta(s)} \mu (ds) }    } \mu (dh).}
 Hence the claim follows
using further the fact that $\mu$ is  translation invariant.  \QED

		\BRM \label{hope}
		Suppose that $Z(t),t\inr^d $ has cadlag sample paths. If for any compact $K\subset \R^d$
		we have that  $\E{\sup_{t \in  K } e^ {Z(t)}} < \IF$, then
		as in the proof above it follows that for any $T>0$
		\BQNY \E[\Big]{ \sup_{t\in [0,T]^d} e^{ Z (t)}}		=      T^d
		\int_{  [0,1]^d} \mathbb E\left\{\frac{\sup_{ \times_{i=1}^d[-h_i T,(1-h_i) T] } e^{  \Theta(t)}}
		{ \int_{   \times_{i=1}^d[-h_i T,(1-h_i) T] } e^{ \Theta(s)}\mu(ds) } \right\}\, \mu^T(dh) ,
		\EQNY
		with  $\mu^T(dh)=\lambda(T dh)/T^d $ where $\lambda$ is  the Lebesgue measure on $\R^d$.}
	\ERM
}

 		\prooftheo{mentalist} 	$i)$ With the same notation as in \eqref{newiPick}, the assumption \eqref{mixing}  implies for any $\delta> 0$ that 
 		\bqny{
 			\limit{T} \frac{\eta_T(h)}{\delta^d}= \mathbb E\left\{\frac{\sup_{t\in \Zdd  } e^{  \Theta (t)}}
 			{ \delta^d \int_{  \R^d  } e^{ \Theta(s)}\mu_\delta(ds) } \right\}< \IF,
 		}
 		where $\mu_\delta$ denotes the counting measure on $\Zdd$.
 		Since  $\delta^d \mu^T (dh)= \delta^d \mu_\delta(T dh)/T^d$ converges weakly as $T \to \IF$ to the Lebesgue measure $\lambda(dh)$, then by \eqref{newiPick}, as in \cite{SBK} we obtain
 		\bqny{
 			\limit{T} \int_{ [0,1]^d}  \eta_T(h) \mu^T(dh)&=&
			 \mathbb E\left\{\frac{\sup_{t\in \Zdd  } e^{\Theta (t)}}
 			{ \delta^d \int_{ \R^d  } e^{ \Theta(s)}\mu_\delta(ds) } \right\}  >0,
 		}
 		hence the first claim follows.  		Next, let $\Gamma_{{k}}(f)= \mathbb{I}\{\inf {\rm argmax}_{ {s}\in  \Zdd \cap [0,T]^d} f({s})= {k}\}$ for some
 		${k} \in \Zdd$. Clearly,
 		$\Gamma_{{k}}(f+c)=\Gamma_{{k}}(f)$ for any constant $c$.  Applying \nelem{lemD}  we obtain (below we set ${u}= (u_1 \ldot u_d)$ and
 		$1- {u}= (1-u_1 \ldot 1 - u_d)$)
 		\BQN
 		\frac1{T^d} \E[\Big]{ \sup_{{t}\in \Zdd \cap [0,T]^d} e^{ Z ({t})}} &= & \frac 1 {T^d} \sum_{{k} \in
 			\Zdd \cap [0,T]^d }
 		\E[\big]{e^{ Z({k})} \Gamma_{{k}}(Z)} \notag\\
 		&=& \frac 1 {T^d} \sum_{{k} \in
 			\Zdd \cap [0,T]^d }
 		\E*{ \Gamma_{{k}}(   L^{{k}} \Theta ) )} \label{irgendw} \\
 		&=& \frac 1 {T^d} \sum_{{k} \in    \Zdd \cap [0,T]^d }
 		\pk*{\sup_{ {s}\in  \Zdd \cap [0,T]^d }  \Theta ({s}-  {k}) =0 }
 		\label{findenw}\\
 		&=& \frac{1}{\delta^d} \int_{  [0,1]^d}
 		\pk*{ \sup_{{t}\in \Zdd \cap [- {u}T, (1-{u})T] }   \Theta ({t})=0 } \delta^d \mu^T(d {u})\notag\\
		&\to & \frac{1}{\delta^d} \int_{[0,1]^d}
 		\pk*{ \sup_{{t}\in \Zdd }  \Theta ({t})\le 0 }   \lambda(d {u}), \quad T \to \IF \label{sLL}\\
 		&= &\frac{1}{\delta^d}
 		\pk*{ \sup_{{t}\in \Zdd }  \Theta({t})= 0 } , \notag
 		\EQN
 		where  \eqref{irgendw} follows by statement $c)$ of \netheo{ThMB} and  \eqref{findenw}  is a consequence of the assumption that $ \Theta$
 		has absolutely continuous fidi's. {Note that by \eqref{mixing}  we have the almost sure convergence
 			$\sup_{\norm{{t}}> T, {t}\in \Zdd}  \Theta({t}) \toas - \IF$ as $T\to \IF$, which implies the convergence in probability
 			$$\sup_{{t}\in \Zdd\cap [- {u}_T T, (1-{u}_T)T]}   \Theta({t}) \toprob \sup_{{t}\in \Zdd }   \Theta({t})\ge
 			 \Theta(\oldemptyset)=0,
 			\quad T\to \IF$$
 			for any ${u}_T$ such that $\limit{T}{u}_T={u}\in (0,1)^d$ and thus
 			$$ \limit{T} \pk*{ \sup_{{t}\in \Zdd\cap [- {u}_TT, (1-{u}_T)T]}   \Theta({t})\le 0 } =
 			\pk*{ \sup_{{t}\in \Zdd}   \Theta({t})\le 0 }.  $$
 			Since further   			$\delta^d \mu^T(dh) $ converges weakly as $T\to \IF$ to $\lambda(dh)$,  then \eqref{sLL} is justified from the validity of \eqref{heidi} below.} \\
  	$ii)$ If $\delta=0$, then by Remark \ref{hope} the proof follows  using further
  	\nelem{intg}. 	\QED

 	\prooflem{tiltGauss} If a.s.\ $Z(h)=0$, then the claim is clear. Suppose therefore that $Z(h)$ has positive variance $\sigma^2(h)>0$.
 	For any distinct $t_1=h,t_2 \ldot t_n \in \TT$ the df of $ {Z} ^\HH =(Z ^\HH (t_1) \ldot Z ^\HH (t_n))$ denoted by $F_h $ is specified by
 	\bqn{\label{talente}
 		F_h(d {x}) &=& F( d {x})  e^{x_1 - \pF},
 		\quad
 		{x}\inr^n,
 	}
 where $F$ is the df of ${Z}=(Z(t_1) \ldot Z(t_n))$. For any ${a}\inr^n$ the df of the rv $({a}, {Z} ^\HH )$ is obtained by tilting the Gaussian rv $({a},  Z)$.  Hence  from here it follows that $ Z^{[h]}$ is Gaussian with the same covariance matrix as $ Z$.
 We calculate next $\E{Z ^\HH (t)}$.  
 	 For any $t\in \TT$
 	\bqny{
 		\E[\big]{ Z ^\HH (t)}&=& \E[\big]{e^{Z(h)- r(h,h)/2} Z(t)  } =\E[\big]{e^{Z(h)- r(h,h)/2} }[r(t,h)+ \E{Z(t)}] =r(t,h)+ \E{Z(t)},
 	}
 	where the second equality  follows by Stein's Lemma, see e.g., (3.4) in \cite{HKume}. 
 	The converse follows easily by \eqref{talente} and is therefore omitted.  \QED


 	\prooftheo{DMElef}  Define $\MZ = \sup_{t\in \TT} e^{Z(t)}, \SZ = \int_{\mathcal{T}} e^{Z(s)} \mu(ds) $, which by our assumption are rv's.  Since $\E{e^{Z(t)}}=1, t\in \TT$, then  $\pk{ \MZ>0}$.  
For $Y_h$ defined in \eqref{yh}, by \eqref{armpush} 
 	\bqn{\label{fFact}
 		\pk{ \mathcal{M}_{Y_h}> 0,\mathcal{S}_{Y_h}=0}= 		\pk{\MZ> 0,\SZ=0} = 0.
 	}
 	Given distinct $t_i\in \mathcal{T}, i \le n$, 
 	using \eqref{fFact} together with  the fact that a.s. $V_h(h)> - \IF$  and $W_h(h)=-\IF$,  for any $ x \inr^n$ we have 
 	\bqn{
 		\lefteqn{- \ln \Hxn }\notag\\
 		&=& \E*{ \frac{ \SZ} {\SZ} \max_{1 \le k \le n} e^{- x_k + Z(t_k)} ;  \MZ> 0 }\notag\\
 		&=& \int_{\TT}  \E*{ \frac{ e^{Z(h)}  } { \SZ } \max_{1 \le k \le n} e^{- x_k + Z(t_k) }
; \mathcal{M}_Z> 0 } \mu (dh)\notag\\
 		&=& \int_{\TT}  \E*{ \frac{ e^{Y_h(h) } }  { \mathcal{S}_{Y_h} }  \max_{1 \le k \le n} e^{- x_k + Y_h(t_k)} ;  \mathcal{M}_{Y_h}> 0 } \mu (dh)\notag\\
 		&=& \int_{\TT}  \E*{ \frac{ e^{Y_h(h) } } { \mathcal{S}_{Y_h} } \max_{1 \le k \le n} e^{- x_k + Y_h(t_k)} ;  \mathcal{M}_{Y_h}> 0,J_h=1 } \mu (dh)\notag \\
&=& \int_{\TT}  \E*{ \frac{ e^{V_h(h) }  } { \int_{\TT}  e^{ V_h(s) } \mu (ds) } \max_{1 \le k \le n} e^{- x_k +  V_h(t_k) } ;
 		  \max_{t\in \TT} e^{V_h(t)  }> 0, J_h=1 } \mu (dh)\notag \\
 		&=& \int_{\TT}  \E*{ \frac{ e^{V_h(h) }  } { \int_{\TT}  e^{ V_h(s)- V_h(h) } \mu (ds) } \max_{1 \le k \le n} e^{- x_k +  V_h(t_k)- V_h(h) } ,   J_h=1 } \mu (dh)
 		\notag\\
 		&=& \int_{\TT}  \E*{ e^{V_h(h) + \ln \pk{ J_h=1}}
 		\max_{1 \le k \le n} e^{- x_k +  V_h(t_k)- V_h(h)   - \ln \bigl( \int_{\TT}  e^{ V_h(s)- V_h(h) } \mu (ds) \bigr) } }\mu (dh) \notag\\
 		&=& \E*{ \int_{\TT}     \max_{1 \le k \le n} e^{- x_k + \ThetaH(t_k) -
 			\ln \bigl(
 			\int_{\TT} e^{  \ThetaH(s)} \mu(ds) \bigr)} \mu (dh)},
\label{clm}
 	}
hence the proof follows. \QED

 	\COM{
 		
 		\prooflem{prop2}
 		For notation simplicity we suppress the subscript $h$ writing simply $J, A, B, V, W$ instead of
 		$J_h, A_h, B_h , V_h$ and $W_h$, respectively. For any $t_1 \ldot t_{n+1} \in \mathcal{T}, t_{n+1}=h$ and $x_1 \ldot x_n \inr$
 		we have with the same arguments as in the proof of \netheo{Th1} (set $p:= \pk{J= 1}$)
 		\BQNY
 		\lefteqn{\pk{\cap_{1 \le j \le n+1}\{ \zeta_{Y}(t_j) \le x_j + \ln k\}\Bigl \lvert \zeta_{Y}(t_{n+1})  > \ln k}}\\
 		&\sim &
 		k [ e^{- }- e^{-}]\\
 		&\sim & \sum_{i=1}^{n+1} e^{-x_i- \ln k} \Psi_i ( )-  \sum_{i=1}^{n} e^{-x_i- \ln k} \Psi_i ( )- \\
 		&\sim & [1- e^{-x_{n+1}}] \Psi_{n+1} ( )
 		\EQNY
 	}

{ 
	\proofkorr{9:20ora} In view of statement $c)$ in \netheo{ThMB}
$ L^h \Xi_\oldemptyset V_\oldemptyset \fdd \Xi_h  V_h$, hence the claim follows by \eqref{DMB2}.
If $\mu$ is the Lebesgue measure on $\TT$, then it follows by \eqref{clm} and 
the translation invariance of $\mu$ that $\eZ$ is stationary. 
\QED
}

 	\prooftheo{minior} $i)$ Exactly as in \cite{DombryE} (therein continuous sample paths are assumed), by Fubini-Tonelli theorem using \eqref{6.6} and \eqref{mua} 
\bqny{
	\int_0^\IF \pk{ u \ThetaH \in A} u^{-2} du  &=& 
	\int_0^\IF  \!\!\! \int_{\TT} f(h)  \mathbb{I}\{f(h) > 0\}  
	\mathbb{I}\{ uf/f(h) \in A\}   \nu( df) u^{-2} du\\
	&=&	\int_{\TT} f(h)  \mathbb{I}\{f(h) > 0\}     \Biggl( \int_0^\IF  
	\mathbb{I}\{ uf/f(h) \in A\}   u^{-2} du \Biggr)   \nu( df)  \\
	&=& \int_{\TT} \mathbb{I}\{f(h) > 0\} \Biggl( \int_0^\IF   
	\mathbb{I}\{ u f \in A\}    u^{-2} du\Biggr)      \nu( df)\\
	&=& \int_{\TT}    \mathbb{I}\{f(h) > 0\}  
	\mathbb{I}\{  f \in A\}  q( df), \quad A \in \mathcal{B}(\mathcal{F}_0),}}
which proves \eqref{eq:min}, hence the claim follows.\\
$ii)$ The proof is the same as that of \cite{MR3024101}[Prop. 4.2].\\
$iii)$ By the assumption 
$$	q_h(A)=\int_0^\IF \pk{  u\ThetaH \in A} u^{-2} du= 
\int_0^\IF \pk{u \mathcal{V} \in A}u^{-2} du$$
for any $A\in \mathcal{B}(\mathcal{F}_0)$,  
hence the claim follows from Remark \eqref{remQ}.\\
\QED

\prooftheo{minor2} 
By the definition of the SPT in \cite{BojanS} (see also \cite{SegersEx17}) and \nelem{prop2}  (see  \eqref{LGa}) we have that 
$\Theta = \Xi_\oldemptyset V_\oldemptyset$ with $V_\oldemptyset \fdd Z \lvert (Z(0)> 0)$. Under this setup it is easy to see that \eqref{pla} is a re-formulation of \eqref{xYx} in terms of SPT. 
\QED

 	\section{Appendix}
 Let $\eZ(t),t\in \TT$ be as in Section 4 where $Z$ has representation \eqref{zH} for some $h\in \TT$, and let $Y$ be a random process given by
\bqn{\label{Yt}Y(t)= J_h A(t)+ (1- J_h) B(t)- \ln \pk{J_h=1}, \quad t\in \mathcal{T},
}
with $A, B, J_h$ being mutually independent  and $ \pk{A(h)=0}= \pk{B(h)=-\IF}=1.$ Denote by $\zeta_Y(t),t\in \TT$ the max-stable process associated to $Y$.
 	
 	\BEL \label{prop2}
 	If $\zeta_Y \fdd   \eZ$, then $A \fdd    \ThetaH$. 	
 	
	\EEL
 	
 	\prooflem{prop2}
 	For notational simplicity we suppress the subscript $h$ writing simply $J$ instead of $J_h$. Let $t_1 \ldot t_{n+1} \in \mathcal{T}, t_{n+1}=h$ be distinct and set 
 	$$c_k= 	\frac{  1  }{ 1- e^{- 1/k}}, \quad  p:= \pk{J= 1}>0, \quad A_j:=A(t_j), \quad  B_j:=B(t_j), \quad K=\{1 \ldot n+1\}. $$
 	With the same arguments as in the proof of \netheo{Th1}, using the fact that
 	$A(h)=0$ and $B(h)=-\IF$ almost surely, we obtain for $k>1$ and  $x_1 \ldot x_n \inr$
 	\BQN \label{LGa}
 	\lefteqn{\pk*{\forall  j \in K : \zeta_{Y}(t_j) \le x_j + \ln k \bigl \lvert \zeta_{Y}(t_h)  > \ln k}}\notag \\
 	&=&
c_k
 	\Biggl(
 	e^{-  \int_{\R}  \bigl[ \pk*{ \exists j\in K  : A_j > x_j + \ln (k/p) -y, J=1} +
 		\pk{ \exists j\in K  : B_j > x_j + \ln (k/p) -y, J=0}   \bigr] e^{-y} \, dy } \notag \\
 	&& -
 	e^{- \frac{1}{kp}\int_{\R} \bigl[
 		\pk{ A_{n+1}> - y , \text{ or }  \exists j\in K  : A_j > x_j -y, J=1}
 		+  \pk{ B_{n+1} > - y, \text{ or }  \exists j\in K : B_j > x_j  -y, J=0} \bigr] e^{-y} \, dy }\Biggr) \notag\\
 	&=& c_k 
 	\Bigl[
 	e^{-  \frac{1}{kp} \int_{\R}  \bigl[ p\pk{ \exists j\in K : A_j > x_j  -y} +
 		(1-p) \pk{ \exists j\in K : B_j > x_j -y}   \bigr] e^{-y} \, dy }  \notag \\
 	&& -
 	e^{- \frac{1}{kp}  \int_{\R} \bigl[
 		p \pk{ 0 > - y  , \text{ or }  \exists j\in K  : A_j > x_j  -y}
 		+  (1-p) \pk{ - \IF > - y  , \text{ or }  \exists j\in K : B_j > x_j  -y} \bigr] e^{-y} \, dy }\Bigr] \notag \\
 	& \to &  \int_{\R}  \Bigl[
 	\pk{y > 0 , \text{ or }  \exists  j\in K  : A_j > x_j  -y} -
 	\pk{ \exists j\in K : A_j > x_j  -y}\Bigr]  e^{-y} \, dy,  \quad k\to \IF \notag \\
 	&= &\pk{ \forall j \in K : A_j+ \mathcal{E}  \le x_j},
 	\EQN
 	with $\mathcal{E}$ a unit exponential rv being independent of $A$, hence since a.s. $A(h)=0$ the proof follows. \QED

Finally, we discuss the asymptotics of $\int_{ \R^d} f_n(x) \nu_n(dx) $ as $n$ tends to infinity.
\BEL \label{intg} Let $\nu_n,n\ge 1$ be positive finite measures on $\R^d,d\ge 1$ which converge weakly as $n\to \IF$ on each set $[-k,k]^d, k\inn$ to some finite measure $\nu$. If $f,f_n,n\ge 1$ is a sequence of  measurable functions on $\R^d$, then for any $k\inn$ we have
\BQN \label{v}
\liminf_{n \to \IF} \int_{ [-k, k] ^d } f_n(x) \nu_n(dx) &\ge&
\int_{ [-k,k]^d} \liminf_{n\to \IF, v\to x} f_n(v) \nu(dx).
\EQN
Assume that for any $u_n\inr^d,n\inn$ such that $\limit{n}u_n=u\in B$ and  $\nu(\R ^d \setminus B)=0$, we have
$\limit{n} f_n(u_n)=f(u)$.  If further $f_n,n\inn$ is uniformly bounded on compacts of $\R^d$, then
\BQN \label{heidi}
\label{losa}\limit{n} \int_{\R^d} f_n(x) \nu_n(dx)= \int_{\R^d} f(x) \nu(dx),
\EQN
provided that
\BQN \label{Par}
\limit{k} \sup_{n\ge 1} \int_{ \R^d \setminus  [-k,k]^d} f_n(x) \nu_n(dx)=0.
\EQN
\EEL
\prooflem{intg} The first claim in \eqref{v} is a special case of \cite{Fatou14}[Th.\ 1.1]. In light of \eqref{Par} the claim in
\eqref{losa} can be established if we show that for any integer $k$
$$\limit{n} \int_{ [-k,k]^d} f_n(x) \nu_n(dx)= \int_{ [-k,k]^d} f(x) \nu(dx)< \IF,$$
which follows directly by \cite{MR2370108}[Lemma 4.2], see also   \cite{Soulier}[Lemma 6.1]. \QED

}
	
\COM
{

 		\def\Z{\mathbb{Z}}
 		\def\xi{\zeta}
 		\def\xiW{\eZ}
 		\section{Appendix}
 		The notion of a mixed moving maxima process on $\R$ defined in
 		\eqref{M3} can be extended to the lattice $\delta \Z$; see for instance
 		Remark 7 in \cite{KabExt}. Suppose that $\{\xi^\delta_ {W}(t), t \in \delta \Z\}$
 		is a stationary max-stable process (with standard Gumbel marginals)
 		given by the construction \eqref{eA} with a process $W$, restricted to $\delta \Z$.
 		Let $F_i^\delta$ be independent copies of a process $F_W^\delta$ on $\delta \Z$
 		with
 		$$\sup_{t\in \delta \Z}  F_W^\delta(t) = F_W^\delta(0) = 0$$
 		almost surely, and
 		\begin{align}\label{Cdelta}
 		C^\delta_W = \left(\E{ \sum_{t\in \delta\Z} \exp( F_W^\delta(t))}\right)^{-1} \in (0,\IF).
 		\end{align}
 		We say that $\xi^\delta_{W} $ admits an M3 representation, if
 		\begin{align*}
 		\xiW ^\delta(t) = \max_{i\geq 0}( F_i^\delta(t - P_i^\delta) + Q_i^\delta), \quad t\in \delta \Z,
 		\end{align*}
 		where $\sum_{i=1}^\infty \epsilon_{(P^\delta_i, Q^\delta_i)}$ is a PPP with intensity
 		$ C^\delta_W \cdot \nu_\delta(dt) \cdot  e^{-x}dx$. Here $\nu_\delta / \delta$ is the
 		counting measure on $\delta \Z$. Below we present the counterpart of Theorem 4.1 in \cite{engelke14}
 		for M3 processes on lattices. \
 		
 		\BEL \label{prop_discr}
 		The process $W^\delta, \delta >0$, the restriction of $W$ to $\delta \Z$,
 		can be expressed in terms of the spectral function $F^\delta$ as
 		\begin{align*}
 		\P(W^\delta \in L) =  C^\delta_W \mathbb{E}\left\{  \sum_{t\in\delta\Z}  \mathbf{1} \left\{ F_W^\delta(\cdot + t) - F^\delta_W(t) \in L\right\} \exp({ F^\delta_W(t)})\right\},
 		\end{align*}
 		which is well-defined probability measure by \eqref{Cdelta}.
 		\EEL
 		\prooflem{prop_discr}
 		The proof goes along the lines of the proof of Theorem 4.1 in \cite{engelke14}. \QED

 		Define in the following
 		$$M_i^\delta:=\sup_{t\in\delta \Z} W_i(t), \quad T_i^\delta:= \inf (\arg\sup_{t\in\delta \Z} W_i(t)),
 		\quad  F_i^\delta(\cdot):= W_i(\cdot+ T_i^\delta)- M_i^\delta, \quad t\in \delta \Z,$$
 		which are well-defined.  Let $\Psi_\delta, \delta \ge 0$ be the PPP
 		\begin{align*}
 		\left\{(T_i^\delta, Q_i +M_i^\delta, F_i^\delta: i\in\N \right\},
 		\end{align*}
 		with intensity measure $ C^\delta_Z \nu_\delta(dt)   e^{-x}dx  \mathbb{P}_{F_Z^\delta}(df) $.

\BT \label{CW}  {For any $\delta>0$ the law  of $F_1^\delta$ is given by
 			$ Z_1 \lvert T_1^\delta =0$
 			and $C_Z^\delta= \frac{1} \delta \pk{ \arg \inf _{t\in \delta \Z} Z(t)=0}.$}
 		\ET
 		
 		\prooftheo{CW}
 		When $Z$ is a Gaussian process, the claims are shown in Theorem 8  and
 		Remark 9 in \cite{KabExt}, respectively. In the general case that $Z$ is not Gaussian,
 		the assumed condition $\E{ \sup_{t \in K} e^{Z(t)}}< \IF$ for any compact subset $K\subset \R$
 		is the additional condition needed to prove the same result as claimed in the aforementioned theorem.
 		In view of \cite{DombryK}, this condition is equivalent with the assumption that $\eZ$ has cadlag sample paths, hence the proof follows.
 		\QED
 		
 	}

 	\section*{Acknowledgments} I am in debt to both referees and the Associate Editor for numerous suggestions/corrections.  In particular, the referees suggested the connection between the change of measure and the extremal functions. Thanks to  K. D\c{e}bicki, T. Dieker, T. Mikosch, I. Molchanov and Ph. Soulier for numerous suggestions.  	Support from  SNSF Grant 200021-166274 is kindly acknowledged.
 	
 	\bibliographystyle{ieeetr}
 	\bibliography{EEEA}
 	\end{document}